\documentclass[12pt]{article}
 \usepackage{amsmath, amsthm, amssymb}
 \usepackage[margin=1in]{geometry}
\newtheorem{theorem}{Theorem}

\newif\ifpdf
	\ifx\pdfoutput\undefined
	\pdffalse 
	\else
	\pdfoutput=1 
	\pdftrue
	\fi

	\ifpdf
	\usepackage[pdftex]{graphicx}
	\else
	\usepackage{graphicx}
	\fi

\begin{document}

\ifpdf
	\DeclareGraphicsExtensions{.pdf, .png, .jpg, .tif}
	\else
	\DeclareGraphicsExtensions{.eps, .jpg}
	\fi

\begin{center}
{\Large {\bf Latent Factor Models for Density Estimation\\}} 
\end{center}
\emph{Suprateek Kundu, Dept. of Biostatistics, UNC Chapel Hill, U.S.A. (skundu@live.unc.edu)}.  \\
\emph{David B. Dunson, Dept. Statistical Science, Duke University, U.S.A. (dunson@stat.duke.edu)}. \\

\vskip 12pt
{\noindent {\bf \quad Abstract:}} Although discrete mixture modeling has formed the backbone of the literature
on Bayesian density estimation, there are some well known disadvantages.  We propose an alternative 
class of priors based on random nonlinear functions of a uniform latent variable with an additive
residual.  The induced prior for the density is shown to have desirable properties including ease of 
centering on an initial guess for the density, large support, posterior consistency and straightforward 
computation via Gibbs sampling.  Some advantages over discrete mixtures, such as Dirichlet process
mixtures of Gaussian kernels, are discussed and illustrated via simulations and an epidemiology 
application. 


\newpage

\noindent{\bf 1.\quad INTRODUCTION}  \\
Nonparametric kernel mixture models are increasingly popular in density estimation and high dimensional data modeling. Kernel mixture models have the form:
\begin{eqnarray}
f(y;G) = \int \mathcal{K}(y;\theta)G(d\theta), \label{eq:NP}
\end{eqnarray}
where $G(\cdot$) is a mixing distribution and $\mathcal{K}(\cdot)$ is a probability kernel. The majority of the nonparametric Bayesian development in this area relies on Dirichlet process (DP) priors (Ferguson, 1973; 1974) for $G$. These models have been generalized to density regression by defining dependence on the covariates $x$ in various ways. M\"uller, Elkanli and West (1996) used a DP mixture of multivariate normals to jointly model the density of the response and predictors to induce a prior on $f(y|x)$. In order to let the parameters of the DP vary over the predictor space $\mathcal{X}$, MacEachern (1999) defined dependent Dirichlet processes (DDP) by assigning stochastic processes on the components in Sethuraman's (1994) DP representation: $G_x=\sum_{i=1}^\infty p_i(x)\delta_{\theta_i(x)}$. De Iorio et al. (2004) proposed a fixed-$p$ DDP, while Griffin and Steel (2006) allowed the weights to depend on predictors.  
Dunson, Pillai and Park (2007) instead used predictor-dependent convex combinations of DP components.

There is also a rich literature on using mixture priors in hierarchical latent variable models. Bush and MacEachern (1996) and Kleinman and Ibrahim (1998) proposed DP mixtures on the distributions of random effects. Fokoue and Titterington (2003) and Fokoue (2005) proposed mixtures of factor analyzers (MFA) corresponding to a finite mixture of multivariate normal kernels with a factor-analytic decomposition of the component-specific covariances. Dunson (2006) used dynamic mixtures of DPs to allow a latent variable distribution to change nonparametrically across groups. More recently, Chen et al. (2009) and Carvalho et al. (2008) proposed nonparametric Bayes MFA allowing an uncertain number of factors. Lee, Lu, and Song (2008) placed a truncated DP on the distribution of the latent variables within a structural equation model (SEM), while Yang and Dunson (2010) proposed a centering approach to ensure identifiability of the latent factor distributions. 

The above approaches have relied on discrete mixture models, which have a number of well known complications motivating alternative methods for modeling unknown densities, such as Polya trees (Mauldin et al, 1992; Lavine, 1992, 1994) and logistic Gaussian processes (LGP) (Lenk 1988, 1991; Tokdar 2007). 
Polya trees have appealing properties in terms of denseness, conjugacy and posterior consistency but have disadvantages in terms of favoring overly spiky densities.  LGP has sound theoretical properties and smoothness of the densities can be controlled through the covariance kernel in the GP.  However, posterior computation is a major hurdle.  Recently, Jara and Hanson (2010) proposed dependent tail-free processes where they modeled the tail-free probabilities with LGP dependent on covariates. Their approach is shown to approximate the Polya tree marginally at each predictor value. An alternative was 
suggested by Tokdar, Zhu and Ghosh (2010) relying on LGP for density regression with dimensionality reduction. 

In this article, we focus on a new approach for nonparametric density estimation and regression that induces a prior on the unknown density through placing a flexible prior on a nonlinear regression function $\theta$ in a latent factor model. The proposed class of models is related to Gaussian process latent variable models (GP-LVM) proposed in the machine learning literature (Lawrence, 2005; Silva and Gramacy, 2010), but our modeling details are different and the focus of this literature has been on nonlinear dimensionality reduction with no consideration of density estimation and associated properties. By using GP priors for $\theta$, we obtain substantial control over the smoothness of the induced densities in a very different manner than that achieved by LGP-based models. Unlike LGP-based models, the proposed model has desirable conjugacy properties facilitating posterior computation. In addition, the method has appealing theoretical properties in terms of large support and posterior consistency.

Relative to some density estimation priors, the proposed latent factor approach is quite
easy to generalize to more challenging settings involving multivariate densities, conditional density estimation,
hierarchical modeling and other complexities.  Although our primary focus in this article is to introduce the
formulation, providing an intuition for how the model works, basic properties and computation, we also
give a flavor of generalizations through a simple conditional density estimation example.  In particular, we
consider a model that induces a prior on the conditional density $f(y|x)$ through joint modeling of the
response and predictors through separate nonparametric latent factor models containing the same latent
variables.  This formulation is completely flexible in the marginal densities, while making strong restrictions
on the dependence to address the curse of dimensionality in a related manner to a copula model.  An
attractive feature of our model is that it naturally allows for incorporation of prior information on the marginal
densities of response and predictors through the mean function of the GP.  The utility of incorporating such
prior information is clear in a reproductive epidemiology application we consider. 
 
\vskip 12pt 

\noindent  {\bf 2.\quad DENSITY ESTIMATION} 

{\noindent {\bf 2.1.\quad Model Specification}} \\
Initially suppose $y_i$ are iid draws from an unknown density $f \in \mathcal{F}$, where $\mathcal{F}$ is the set of densities on $\Re$ with respect to Lesbesgue measure.  We propose to induce a prior $f \sim \Pi$ through 
\begin{eqnarray}
y_i & = & \mu(x_i) + \epsilon_i,\quad \epsilon_i \sim \Gamma_\sigma, \nonumber \\
\mu  \sim  \Pi^* ,\quad \sigma & \sim & \nu,\quad x_i  \sim  \mbox{Uniform}(0,1), \label{eq:base}
\end{eqnarray}
where $\mu \in \Theta$ is an unknown $\left[0,1\right] \to \Re$ function, $x_i$ is a uniformly distributed latent variable, and the error distribution $\Gamma_\sigma$ is centered at $0$ and has scale parameter $\sigma$. Hence, in the special case in which $\mu(x) = \mu$, so that the regression function is a constant, and $\Gamma_\sigma$ is normal, we have $f(y; \mu, \sigma^2) = N(y; \mu, \sigma^2)$ so we obtain a normal density. The density of $y$ conditionally on the unknown regression function $\mu$ and $\sigma$ is obtained on marginalizing out the latent variable as 
\begin{eqnarray}
f(y; \mu, \sigma) = f_{\mu,\sigma}(y)=\int_0^1 \Gamma_\sigma(y-\mu(x)) dx. \label{eq:cGPT}
\end{eqnarray}
To complete the specification, we let $\mu \sim \Pi^*,\mbox{ } \sigma\sim \nu$ and obtain the marginal density 
\begin{eqnarray}
f(y) = \int_0^\infty \int_{\Theta} \int_0^1\Gamma_\sigma(y-\mu(x)) dx\Pi^*(d\mu)\nu(d\sigma).  \label{eq:GPT}
\end{eqnarray}
Hence, a prior $f \sim \Pi$ is induced through assigning independent priors to $\mu$ and $\sigma$ in expression (\ref{eq:cGPT}).  When the prior on $\mu$ is a Gaussian process and the error distribution is $N(0,\sigma^2)$ (denoted as $\Gamma_{\sigma}=\phi_\sigma$), we refer to $f_{\mu,\sigma}$ as a Gaussian process transfer (GPT) model and the induced prior $f \sim \Pi$ as a GPT prior.  

The GPT prior does not have the kernel mixture form (\ref{eq:NP}).  There will be no clustering of subjects or label switching issues.  Instead, the prior $f \sim \Pi$ is induced through adding a Gaussian residual to a Gaussian process regression model in a uniform latent variable.  This is a simple structure aiding computation and interpretability.  One can control the smoothness of the density through the covariance in the GP prior for the regression function $\mu$ and the size of the scale parameter $\sigma$.  In limiting cases, one can obtain realizations of $\mu$ concentrated close to a flat line, leading to a normal density as a special case.  In addition, by making $\sigma$ small and choosing the GP covariance to generate a very bumpy $\mu$, one can obtain arbitrarily bumpy densities.  In practice, by choosing hyperpriors for key covariance parameters, we obtain a data adaptive approach that often outperforms discrete kernel mixtures.  The performance of discrete kernel mixtures relies on the ability to accurately approximate the density with few components, and DP mixtures tend to heavily favor a small number of dominate kernels.  This tendency can sometimes lead to relatively poor estimation, as illustrated in section 6.
\vskip 12pt 

\noindent {\bf 2.2.\quad Prior Specification}  \\
Prior elicitation is an important aspect of Bayesian modeling, with the prior playing a particularly important role in Bayesian nonparametric models involving infinitely many parameters.  Most of the Bayesian nonparametrics literature relies on default priors, which do not reflect available prior knowledge in a particular application area, but are chosen to lead to good performance in terms of posterior behavior in a wide variety of applications.  However, as in parametric models, well chosen informative priors that utilize information, such as historical data on the variables under study, can substantially improve the performance in small to moderate samples.  In DP mixtures, such prior information is typically incorporated through choice of hyperparameters in the base measure, while maintaining conjugacy for ease in computation.  For example, in Gaussian kernel mixtures, a normal-inverse gamma base measure would be chosen having parameters representing prior knowledge.  This is appropriate when prior knowledge implies that the density follows a $t$ distribution, but when one has prior information that the density follows a more complex form (as in our premature delivery application) then elicitation is substantially more difficult.  Obtaining a base measure that leads to a particular elicited density is a deconvolution problem, which can be difficult to solve for non-atomic base measures.  In addition, posterior computation under the resulting complex and non-conjugate base measure may be challenging.  An advantage of the GPT is that the prior for the density can be centered on an arbitrary choice easily through the prior mean in the GP prior for $\mu$.

To elaborate, Theorem 3 (section 3) ensures that $f_{\mu,\sigma}\approx f_{\mu_*,\sigma}$, when $\mu \approx \mu_* = F_*^{-1}$ and $\sigma \approx 0$. In terms of application, this translates to incorporating a prior guess $f_{\mu_*,\sigma}$ for the density through the corresponding mean function $\mu_*=F_*^{-1}$ of the GP, and letting the prior for $\sigma$ to have mode near 0. Such a mean function can be constructed by obtaining frequentist kernel estimates of the concerned density using some external data, and then converting it into an inverse cdf on a grid of points in [0,1] (using a linear approximation). Thus, the characteristics of the entire density is captured through $F_*^{-1}$ (mean function of the GP) and we let the data influence the deviation of the posterior from the prior guess. 

These ideas are demonstrated in Figure 1, where we use some earlier data on gestational age at delivery to construct prior densities. We choose a Ga(25,1) prior for the residual precision, and different sets of hyper-parameters for the covariance kernel of the GP. The frequentist kernel estimates were obtained by the bandwidth selection method of Sheather and Jones (1991), using a Gaussian kernel (`kernel' function in R). It is evident that the smoothness as well as the degree of deviation of the prior from the frequentist estimate can be controlled through the hyper-parameters in the covariance kernel of the GP. 

\noindent {\bf 3.\quad Theoretical Properties}  \\
To further justify the proposed prior, we show large support and posterior consistency properties.  Large support is an important property in that it ensures that our prior can generate densities that are arbitrarily close to any true density $f_0$ in a large class, a defining property for a nonparametric Bayesian procedure and a necessary condition to allow the posterior to concentrate in small neighborhoods of the truth.  Instead of focusing narrowly on GPT priors, we provide broad theoretical results for priors in the general class of expression (\ref{eq:base}).  

Before proceeding, it is necessary to define some notation and concepts.  We denote the Kullback-Leibler (KL) divergence of $f_{\mu,\sigma}$ from $f_0$ as $KL(f_{\mu,\sigma},f_{0})$ and an $\epsilon-$sized KL neighborhood around $f_0$ as $KL_{\epsilon}(f_0)$. The sup-norm distance is denoted by $||.||_\infty$. Our development will rely on the fact that any density $f_0$ can be calculated from a ``true function'' $\mu_0 = F_0^{-1}$, with $F_0$ denoting the cumulative distribution function.  To generate $y_i \sim f_0$, one can equivalently draw $x_i \sim \mbox{Uniform}(0,1)$ and let $y_i = \mu_0(x_i)$.  This corresponds to the limiting case as $\sigma \rightarrow 0$ in model (\ref{eq:base}) with $\mu = \mu_0$. Hence, we assume a weak regularity condition that the true density can be represented as the limiting case
\begin{equation}
f_{0}\left(y\right)= \lim_{\sigma \rightarrow 0} \int_{0}^{1} \Gamma_{\sigma}(y-\mu_{0}(x)) dx ,\label{eq:true}
\end{equation}
assuming $\Gamma_{\sigma}$ is chosen so that such a limit exists.  The above condition is quite reasonable, only excluding densities for which convergence in distribution does not imply convergence of the corresponding density functions. As additional reasonable regularity conditions, we assume that $f_0$ is strictly positive and finite, sup$_x|\mu_0(x)|<\infty$ and $0<\Gamma_\sigma(u)<\infty$ for finite $u$. A uniformly bounded $\mu_0$ along with the condition on $\Gamma_\sigma$ implies that for $\mu$ belonging to a ball of finite radius around $\mu_0$, $f_{\mu,\sigma}$ is strictly positive and finite for all $\sigma \in \Re^+$, which ensures a finite KL divergence for a suitable subset of $\mu$ values.

\begin{theorem} Let sup$_x|\mu_0(x)|<\infty$, $\Gamma_{\sigma}$ be normal, Laplace or Cauchy with scale parameter $\sigma$ and $f_0$ be the corresponding density in $\mathcal{F}$ defined as in equation ($\ref{eq:true}$). If $\mu_0$ is in the sup-norm support of $\Pi^*$ and \small{$\mbox{ }\nu\Big\{ \sigma:\sigma \in (0,\eta)\Big\}>0$} for all $\eta>0$, then $\Pi(KL_\epsilon(f_0))>0$ for all $\epsilon>0$. 
\end{theorem}

Theorem 1 allows us to verify that the induced prior on the density $f$ assigns positive probability to KL neighborhoods of any strictly positive and finite true density $f_0$.  From Schwartz (1965), if the true density $f_0$ is in the KL support of the prior for $f$, the posterior distribution for $f$ will concentrate asymptotically in arbitrarily small weak neighborhoods of $f_0$. Theorem 1 requires the prior $\mu \sim \Pi^*$ to place positive probability in sup-norm neighborhoods of the inverse cdf $F_0^{-1}$.  Although one can verify this condition for certain choices of $\Pi^*$, such as appropriately chosen Gaussian process priors, it is nonetheless somewhat stringent.  We show in Theorem 2 that this condition can be relaxed to only require that the prior $\mu \sim \Pi^*$ assigns positive probability to L-1 neighborhoods of any element $\mu_0$ of $\Theta$.  It is well known that positive sup-norm support automatically guarantees positive L-1 support but the converse is not true.

\begin{theorem} Let sup$_x|\mu_0(x)|<\infty$, $\Gamma_{\sigma}=\phi_\sigma$ and f$_{0}$ be the corresponding density in $\mathcal{F}$ defined in equation ($\ref{eq:true}$). If $\mu_0$ is in the L-1 support of $\Pi^*$ and \small{$\mbox{ }\nu\Big\{ \sigma:\sigma \in (0,\eta)\Big\}>0$} for all $\eta>0$, then $\Pi(KL_\epsilon(f_0))>0$ for all $\epsilon>0$. 
\end{theorem}

As the prior $f \sim \Pi$ is specified indirectly through priors $\mu \sim \Pi^*$ and $\sigma \sim \nu$, it is desirable for elicitation purposes to verify that, for sufficiently small $\sigma$, $\mu \approx \mu_0$ implies that $f_{\mu,\sigma} \approx f_0$.  Theorem 3 provides such a verification assuming Gaussian errors.  This implies one can potentially center the prior for the density $f$ on an initial parametric guess $\tilde{f}$ by centering $\mu \sim \Pi^*$ on the inverse cdf $\tilde{F}^{-1}$ while choosing the prior for $\sigma$ to have mode near zero.  The data will then inform about the degree to which $\mu$ deviates from $\tilde{F}^{-1}$ and $\sigma$ deviates from 0.  

\begin{theorem} For $\mu_0 \in \Theta$ and $\Gamma_{\sigma}=\phi_\sigma$, let $f_0$ be the density resulting from equation ($\ref{eq:true}$). Then for $\mu \in N_{\epsilon_{1}}\left(\mu_0\right)$, with $N_{\epsilon_1}(\mu_0)$ an $\epsilon_1$-sized L-1 neighborhood around $\mu_0$, and $\sigma \in (\epsilon_{2},\mbox{ } \epsilon_{2}^{*})$,  we have $f_{\mu,\sigma} \in N_{\frac{\epsilon_{1}}{\epsilon_{2}}}\left(f_{0}\right)$ for arbitrarily small $\epsilon_1,\epsilon_2,\epsilon_2^*$ such that  $0<\epsilon_{1}<\epsilon_{2}<\epsilon_{2}^{*}$.
\end{theorem}

Although Theorems 1-2 lead to weak posterior consistency, small weak neighborhoods around $f_0$ are topologically too large and may include densities that are quite different from $f_0$ in shape and other characteristics.  Hence, it is appealing to establish a strong posterior consistency result in which the posterior probability allocated to arbitrarily small L-1 neighborhoods of $f_0$ increases towards one exponentially fast with increasing sample size.  Focusing on the GPT prior described above, we show in Theorem 4 that strong posterior consistency holds under some conditions on the prior. Notably, for an appropriately chosen GPT prior, we obtain L-1 posterior consistency for all strictly positive and finite true densities $f_0 \in \mathcal{F}$ with the weak regularity condition ($5$). 

\begin{theorem} Suppose $\mu\sim GP(m,c)$ and define $f_0$ as in equation ($\ref{eq:true}$) where $\Gamma_{\sigma}=\phi_\sigma$. Let U=$\left\{f_{\mu,\sigma}:\int|f_{\mu,\sigma}-f_0|dy<\epsilon,\mu \in \Theta, \sigma \in (0,\infty) \right\}$. Suppose the mean function m($\cdot$) is continuously differentiable, and the covariance function c($\cdot,\cdot$) has continuous fourth derivatives. Further $0<\nu(\sigma:\sigma \in (0,\eta_1))<\eta_2$ for any arbitrarily small $\eta_1$, with $\eta_2$ decreasing with $\eta_1$. Then, $f_0$ is in the KL support of $\Pi$ implies that the posterior is strongly consistent at $f_0$.
\end{theorem}
The above assumptions on the GP can be verified for many popular covariance functions, both stationary and nonstationary. Some such examples can be found in Choi et. al. (2004).
\vskip 12pt 

{\noindent {\bf 4. \quad SINGLE FACTOR DENSITY REGRESSION}} \\
As a simple and parsimonious single factor model that generalizes the model of Section 2 to include predictors $z_i = (z_{i1},\ldots, z_{ip})'$ of a response $y_i$, we let 
\begin{eqnarray}
y_i & = & \mu^Y(x_i) + \epsilon_i,\quad \epsilon_i \sim N(0,\sigma_Y^2), \nonumber \\
z_{ik} &=&  \mu^{Z_k}(x_i) + \epsilon^*_{ik},\quad \epsilon^*_{ik} \sim N(0,\sigma_{Z_k}^2), k=1,2,\ldots,p, \nonumber \\
x_i & \sim &  \mbox{Uniform}(0,1), \nonumber \\
\mu^Y & \sim & \Pi^Y , \quad \mu^{Z_k} \sim \Pi^{Z_k}, k=1,2,\ldots,p, \nonumber \\
\sigma^Y & \sim & \nu,\quad \sigma_{Z_k} \sim \nu, \label{eq:regression}
\end{eqnarray}
where $\mu^Y,\mu^{Z_k} \in \Theta$ are unknown $\left[0,1\right] \to \Re$ functions, $\epsilon$'s are independent errors and x$_i$ is the latent variable. For simplicity, we assume the same prior $\nu$ on the precision of the measurement errors in each component model, though this assumption is trivial to relax.  Expression (\ref{eq:regression}) is a multivariate generalization of the univariate density estimation model (\ref{eq:base}).  Marginally each of the variables is assigned exactly the prior in (\ref{eq:base}) and to allow dependence we incorporate the same latent factor $x_i$ in each of the models.   

Our goal in defining a joint model is to induce a flexible but parsimonious model for the conditional density of $y_i$ given the predictors ${\bf z}_i$.  In estimating conditional densities for multiple predictors, one encounters a daunting dimensionality problem in that one is attempting to estimate a density nonparametrically while allowing arbitrary changes in this density across a multivariate predictor space.  Clearly, as $p$ increases even for large samples there will be many regions of the predictor space that have sparse observations.  As a compromise between flexibility and parsimony in addressing the curse of dimensionality, we propose to use a single factor model in which the marginals for each variable are fully flexible but restrictions come in through assuming dependence on a single $x_i$.  Extensions to the multiple factor case are straightforward.  

\vskip 12pt 

\noindent {\bf 5. \quad POSTERIOR COMPUTATION}  \\
For simplicity, we focus on the single predictor density regression case when outlining an MCMC algorithm for posterior computation.  
Let Y$_{n\times1}$ and Z$_{N\times1}$ denote the vector of observations and covariates, respectively.  We are interested in prediction of $y_{n+1},\ldots,y_N$ based on $z_{n+1},\ldots,z_N$. Let $\mu^n_Y$ ($n \times 1$) and $\mu^N_Z$ ($N \times 1$)  denote the realizations of the GP $\mu^Y$ and $\mu^Z$ at the latent variable values ${\bf x} = (x_1,\ldots, x_n,x_{n+1},x_N)'$. From the GP prior, we have $\mu^n_Y \sim \mbox{N}_n( m_Y^n, \textbf{K}_Y^n)$ and $\mu^N_Z \sim \mbox{N}_N( m_Z^N, \textbf{K}_Z^N)$. The covariance kernels are squared exponential with \textbf{K}$_Y(x,x')=\frac{1}{\phi_Y}\exp\bigg\{-C_Y(x-x')^2\bigg\}$ and \textbf{K}$_Z(x,x')=\frac{1}{\phi_Z}\exp\bigg\{-C_Z(x-x')^2\bigg\}$. We specify conjugate gamma priors:  $\sigma_Y^{-2}\sim Ga(a_\sigma,b_\sigma)$, $\sigma_Z^{-2}\sim Ga(aa_\sigma,bb_\sigma)$, $\phi_Y \sim Ga(a_\phi,b_\phi)$ and $\phi_Z \sim Ga(aa_\phi,bb_\phi)$. For updating the latent variables {\bf x}, we adopt the griddy Gibbs approach using a set of 
 evenly distributed grid points $g^*_1,g^*_2,\ldots,g^*_G \in (0,1)$. Let \textbf{D}$_Y$ and \textbf{D}$_Z$ be diagonal matrices having $\sigma^2_Y$ and $\sigma^2_Z$ as their diagonal elements respectively. Let $\mu^n_Y(-i)$ include all elements of $\mu^n_Y$ except $\mu^Y(x_i)$, and similarly for $\mu^N_Z(-i)$. The Gibbs sampling algorithm alternates between the following steps.\\
\emph{Step1}: Update $\sigma^2_Y$ and $\sigma^2_Z$ using 
\begin{small}
$\pi(\sigma^{-2}_Y|-)\sim$Ga(a$_\sigma$+n/2, b$_\sigma + \frac{1}{2}\sum_{i=1}^n(y_i- \mu^Y(x_i))^2$) and $\pi(\sigma^{-2}_Z|-)\sim$Ga(aa$_\sigma$+N/2, bb$_\sigma + \frac{1}{2}\sum_{i=1}^N(z_i- \mu^Z(x_i))^2$)
\end{small} 
respectively.\\ 
\emph{Step2}: To sample the latent variables, choose $x_i=g_k^*$ with probability p$_{ik}$, where 
\begin{eqnarray*}
p_{ik}=P[x_i=g_k^*|-]&=& 
\frac{p_{ik}^Y p_{ik}^Z N(\mu^Y(x_i=g_k^*)|\mu^n_Y(-i))  N(\mu^Z(x_i=g_k^*)|\mu^N_Z(-i))}{\sum_{l=1}^G p_{il}^Y p_{il}^Z  N(\mu^Y(x_i=g_l^*)|\mu^n_Y(-i)) N(\mu^Z(x_i=g_l^*)|\mu^N_Z(-i))}, \mbox{ if i}\le \mbox{ n}\\
&=& \frac{p_{ik}^Z  N(\mu^Z(x_i=g_k^*)|\mu^N_Z(-i))}{\sum_{l=1}^G p_{il}^Z  N(\mu^Z(x_i=g_l^*)|\mu^N_Z(-i)) }, \mbox{ if n}<i\le \mbox{ N},
\end{eqnarray*}
where $p_{ik}^Y=N(y_i;\mu^Y(x_i=g_k^*), \sigma^2_Y)$, $p_{ik}^Z=N(z_i;\mu^Z(x_i=g_k^*), \sigma^2_Z)$ and k=$1,2,\ldots,G$.\\
\emph{Step3}: Update $\mu^n_Y$ and $\mu^N_Z$ using 
\begin{small}
$\pi(\mu^n_Y|-)=N_n\Big( (\textbf{D}_Y^{-1} + (\textbf{K}_Y^n)^{-1})^{-1}(\textbf{D}_Y^{-1}Y + (\textbf{K}_Y^n)^{-1}m^n_Y), (\textbf{D}_Y^{-1} + (\textbf{K}_Y^n)^{-1})^{-1}\Big)$ and $\pi(\mu^N_Z|-)=N_N\Big( (\textbf{D}_Z^{-1} + (\textbf{K}_Z^N)^{-1})^{-1}(\textbf{D}_Z^{-1}Z + (\textbf{K}_Z^N)^{-1}m^N_Z), (\textbf{D}_Z^{-1} + (\textbf{K}_Z^N)^{-1})^{-1}\Big)$
\end{small}
respectively.\\
\emph{Step4}: Update $\mu^{*G}_Y=\left\{\mu^Y(g^*_1),\ldots,\mu^Y(g^*_G)\right\}$ and $\mu^{*G}_Z=\left\{\mu^Z(g^*_1),\ldots,\mu^Z(g^*_G)\right\}$ using the conditional normal distributions N($\mu^{*G}_Y|\mu^n_Y$) and N($\mu^{*G}_Z|\mu^N_Z$) respectively. \\
\emph{Step5}: Update $\phi_Y$ and $\phi_Z$ using $\pi(\phi_Y|-)\sim$Ga(a$_\phi+\frac{n}{2}$, b$_\phi + \frac{1}{2}(\mu^n_Y-m^n_Y)'(\textbf{K}_Y^n)^{-1}(\mu^n_Y-m^n_Y)$ ) and $\pi(\phi_Z|-)\sim$Ga(aa$_\phi+\frac{N}{2}$, bb$_\phi + \frac{1}{2}(\mu^N_Z-m^N_Z)'(\textbf{K}_Z^N)^{-1}(\mu^N_Z-m^N_Z)$ ) respectively. \\ 
\emph{Step6}: Update C$_Y$ and C$_Z$ using Metropolis random walk for log(C$_Y$) and log(C$_Z$). \\
For prediction of $y_k$ based on $z_k$, k = $n+1,\ldots,N$, we use $\pi(y_k|-)=N(y_k;\mu^Y(x_k),\sigma^2_Y)$, while the conditional density estimate is calculated as $ \hat{f}(y|z)= \frac{\frac{1}{G}\sum_{k=1}^G \phi_{\sigma_Y}(y-\mu^Y(g^*_k) )\phi_{\sigma_Z}(z-\mu^Z(g^*_k) )}{\frac{1}{G}\sum_{k=1}^G\phi_{\sigma_Z}(z-\mu^Z(g^*_k) )}$. 

\vskip 12pt 

\noindent  {\bf 6. \quad SIMULATION STUDY}  \\
To assess the performance of the GPT approach in density estimation as well as density regression, we conducted several simulation studies. We chose the mean function for the GP as m(x)=2sin(x)+cos(x) and utilized the squared exponential covariance kernel. For computational purposes, we worked with the standardized data and then transformed it back in the final step. The hyperparameters for the gamma priors were chosen to be one throughout. Although we used 75 grid points for the griddy Gibbs approach, the number of points could be as low as 60. The number of iterations used was 10000 with a burn in of 1000. The convergence for the main quantities such as $\mu$ was rapid with good mixing. All results are reported over 5 replicates. \\
\noindent {\bf 6.1.\quad Univariate Density Estimation}  \\
To see how well the GPT does in practice for density estimation, we looked at a variety of scenarios, where the truth was generated from the densities considered in Marron and Wand (1992), which are essentially finite mixtures of Gaussians. We present the results from four of those cases which we thought to be interesting deviations from normality and could be potentially encountered in applications. These are the 2nd, 6th, 8th and 9th Marron-Wand densities. The sample size used was 100. For comparison, we looked at DP mixture of Gaussians (Escobar and West, 1995), mixtures of Polya trees (Hanson, 2006) and frequentist kernel estimates using a Gaussian kernel (and the bandwidth selection method of Sheather and Jones, 1991). More specifically, for both DP mixtures and mixtures of Polya trees, we used the DP package in R and the standard hyperparameter values therein. We used algorithm 8 of Neal (2000) with m=1 for DP mixtures of Gaussians. For frequentist kernel, we used the function ``density" in R with Gaussian kernel. Overall, we found that varying the hyperparameter values within a reasonable range does not significantly alter the density estimation results for a sample size of 100, for any of the competitors. Table 1 presents the L-1 distance between true and estimated densities while Figure 2 depicts the density plots.\\
\begin{small}
\textbf{Table 1: Marron-Wand Curves: L-1 distance between true and estimated densities}
\begin{center}
\begin{tabular}{|c c c c c|}
\hline 
 Method & &L-1 Distance & & \\
 &MW 2 & MW 6 &MW 8 &MW 9 \\[0.08in]

   GPT     \quad &0.031 \quad &0.035 \quad &0.031 \quad &0.028 \\
                                        
   DPM     \quad &0.035 \quad &0.036 \quad &0.03 \quad &0.038  \\
                                
   Polya tree mixture    \quad &0.065 \quad &0.036 \quad &0.045 \quad &0.042 \\
                               
   Frequentist Kernel  \quad &0.145 \quad &0.031 \quad &0.033 \quad &0.028  \\  
 
\hline
\end{tabular}
\end{center}
\end{small}
From table 1, we see that even when the truth is generated from a finite mixture of Gaussians, the GPT tends to do better or at least as well as the DP mixture of Gaussians. Mixtures of Polya trees have somewhat worse performance and result in overly spiky looking estimates. \\
\noindent {\bf 6.2.\quad Single Factor Density Regression}  \\
For density regression, we generated a univariate response by allowing the conditional mean as well as the residual error distribution to vary with the covariate. We compared the out of sample predictive performance of GPT with other competitors such as DP mixture of bivariate normals (M\"uller, Erkanli and West, 1996), Bayesian additive regression trees (BART) (Chipman, George and McCulloch, 2010), GP mean regression (O'Hagan and Kingman, 1978) and treed GP (Gramacy and Lee, 2008), based on standard packages in R. We used the DP package for DP mixtures of Gaussians and the Bayestree package for the other three methods, and the hyperparameter values therein. The density regression results did not change significantly on varying the hyperparameter values within a reasonable range, for all the competitors. We used the following scheme for simulations:
\begin{eqnarray*}
 Z \sim F_Z,  \quad y_{i} &=& \lambda \exp\bigg(-\frac{e^{z_i}}{1+e^{z_i}}\bigg) + \frac{e^{z_i}}{1+e^{z_i}}\epsilon_{i}, \quad \epsilon_{i} \sim N(0,\sigma^2),
\end{eqnarray*}
where $F_Z$ is the distribution of the predictors which was chosen to be a trimodal density (9th Marron-Wand curve). We chose $\lambda=3$ and split the total sample size of 100 into training set of 50 and test set of 50. The above data generating model allows the shape of the conditional density to change with predictors, hence making prediction non-trivial. Table 2 shows the performance of the GPT along with a few competitors. We computed the mean square error (MSE), 95\% coverage for the mean (COV), as well as the L-1 distance between true and estimated densities at 25th, 50th and 75th percentiles of the predictor distribution. \\
\begin{small}
\textbf{Table 2: Out of Sample MSE and L-1 distance between true and estimated densities}
\begin{center}
\begin{tabular}{|c c c c c c |}
\hline
 
 Method &MSE &COV(\%)  & &L-1 Distance &  \\
\hline
  & & &25th &50th &75th \\[0.08in]

   GPT     \quad &1.26 \quad &94  \quad &0.08 \quad &0.04 \quad &0.06  \\
                                        
   DPM     \quad &1.53 \quad &42  \quad &0.03 \quad &0.04 \quad &0.06  \\
                                
   BART    \quad &1.59 \quad &46  \quad &0.13 \quad &0.026 \quad &0.10  \\
                               
   GP reg  \quad &1.52 \quad &76  \quad &0.14 \quad &0.03 \quad &0.10  \\
                               
   treed GP \quad &1.6 \quad &72  \quad &0.07 \quad &0.09 \quad &0.04  \\                             
 
\hline
\end{tabular}
\end{center}
\end{small}

The results in table 2 are consistent with our experience in simulations- when the predictor distribution is multimodal and the shape of the conditional density is allowed to change with predictors, then the GPT tends to do as well or better than DP mixture of Gaussians. For the above study, the average number of components in the conditional distribution obtained from DP mixtures was around 15 which is quite high for a sample size of 50. As illustrated in table 2, BART, treed GP and the GP mean regression methods are primarily mean regression methods and so cannot possibly do well in terms of characterizing the entire conditional of response given predictors. They might perhaps estimate the mean surface reasonably well, but eventually fail in capturing multimodality or tail behavior, the latter often being an important focus of inferences.
\vskip 12pt 

\noindent  {\bf 7. \quad EPIDEMIOLOGY APPLICATION} \\
{\noindent {\bf 7.1 \quad Study Background}}\\
DDT is a cheap and popular alternative for reducing the transmission of malaria, but has been shown to have negative effects on public health. In order to study the association between the DDT metabolite DDE and preterm delivery, Longnecker et al. (2001) measured DDE in mother's serum in the third trimester of pregnancy and also recorded gestational age at delivery (GAD) as well as age. They did logistic regression with response as dichotomized GAD (preterm or normal depending on a cut-off of 37 weeks of completed gestation) and explanatory variables as categorized DDE based on empirical quantiles. Their results showed a significant dose-response relationship which had important public health implications. Dunson et al. (2008) analyzed the data using kernel stick-breaking processes, and showed an increasing bump in the left tail of the GAD density with increasing DDE. 

{\noindent {\bf 7.2 \quad Analysis and Results} \\
We used the GPT to analyze the dose response relationship in a subset of 182 women of advanced maternal age ($\ge 35$ yrs) in the above dataset. We examined the conditional distribution of GAD at 10th, 60th, 90th and 99th percentile of DDE. Further, we looked at the dose response relationship between preterm birth and DDE, by examining the left tail of GAD over varying doses of DDE. We used normalized data for analysis and converted it back in the final step. Using the prior specification approach of section 2.2, we were able to incorporate prior information on the marginal density of GAD (using an external data) through the mean function of the GP. Note that prior on $\sigma^{-2}$ for GAD was chosen as Ga(25,1). Given the limited sample size and the complexity of the data we are trying to model, we adjusted other hyperparameter settings to reflect our prior belief about the data. The starting value for the length-scale parameter in the covariance kernel in the Metropolis random walk was chosen to be 25, so as to have smooth Gaussian process prior. Instead of working with DDE, we used log(DDE) which resembled a Gaussian distribution, with a 0 mean function for the predictor component and Ga(1,1) prior for the corresponding residual precision.

Figure 3 shows the conditional distribution curves for GPT along with 90\% credible intervals. Although we focused on a small subsample of 182 women of advanced maternal age, the GPT results for the conditional density are remarkably similar to the ones reported in Dunson et al. (2008), which suggests that there is no systematic difference for women of advanced maternal age. The conditional densities show an increasing bump in the left tail with increasing DDE, suggesting increased risk of preterm birth at higher doses. This is further supported by dose-response curves for P(GAD$<$T) in Figure 5, with different choices for cut-off T. Although the dose-response curve is mostly flat for T=33 weeks, the relationship becomes more significant as cut-off increases, with the dose-response tapering off at T=40 weeks. This suggests that increased risk of preterm birth at higher DDE dosage is attributable to premature deliveries between 33 and 37 weeks. Trace plots of f($y|z$) for different DDE percentiles (not shown) exhibit excellent rates of convergence and mixing. For comparison, Figure 4 shows the density estimates from the DP mixture of Gaussians which has a tendency to overly favor multimodal densities, which is as expected given our simulation study results. These results were obtained using DP package in R (and the data driven hyperparameter values therein), which utilizes algorithm 8 of Neal (2000) with m=1.  \\
\vskip 12pt 
\noindent  {\bf 8. \quad Discussion}  \\
In this paper, we propose a latent factor model for density estimation. This novel method provides us with a flexible non-discrete mixture alternative to be used in a variety of situations including density estimation, density regression, hierarchical latent variable models and even mixed models. We provide theoretical theoretical justifications for GPT and demonstrate it's usefulness as a building block for more complex models involving covariates. Building on our work, Pati, Bhattacharya and Dunson (to be submitted) recently showed minimax optimal rates of posterior contraction for Bayesian density estimation from non-linear latent variable models, also obtaining initial results on contraction rates in conditional density estimation.  The close relationship between non-linear latent variable models for densities and non-linear mean regression models facilitates not only posterior computation but also derivations of theoretical properties, such as contraction rates, which have proven difficult to study for discrete mixtures beyond simple settings. 
\vskip 12pt 
{\noindent  {\bf 9. \quad Acknowledgments} } \\
The authors thank Debdeep Pati and Anirban Bhattacharya for their helpful comments.  This work was support by Award Number R01ES017240 from the National Institute of Environmental Health Sciences.  The content is solely the responsibility of the authors and does not necessarily represent the official views of the National Institute of Environmental Health Sciences or the
National Institutes of Health.
\vskip 12pt
{\noindent {\quad APPENDIX: PROOF OF RESULTS}}\\
\textbf{Proof of Theorem 1}:\\
That f$_{\mu,\sigma}$ exists is seen from equation ($\ref{eq:cGPT}$). The limiting case as $\sigma \rightarrow 0$ can be explained as the case when Y has the distribution function F$=\mu^{-1}$ and marginalizing out the latent variable x$\in [0,1]$. Thus 
\begin{small}
$\lim_{\sigma \rightarrow 0} \int_{0}^{1} \Gamma_\sigma(y-\mu(x)) dx $ 
\end{small}
exists. Under regularity conditions, we can write the KL divergence between f$_{\mu,\sigma}$ and f$_0$ as 
\begin{eqnarray}
KL(f_{\mu,\sigma}, f_0)&=&\int f_{0} \log\frac{f_0}{f_{\mu,\sigma}}=\int f_0(y) \left[\frac{\lim_{\sigma \rightarrow 0}\int_{0}^{1} \Gamma_{\sigma}(y-\mu_{0}(x)) dx}{\int_{0}^{1} \Gamma_{\sigma}(y-\mu(x)) dx }\right]dy \label{eq:weak}. 
\end{eqnarray}
Note,  $\Gamma_{\sigma}\left(y-\mu(x)\right)\ge\frac{\Gamma_{\sigma}\left(y-\mu_0(x)\right)}{h_\sigma\left(y,\mu(x)-\mu_0(x)\right)}$, where \begin{small}
h$_\sigma(y,\mu-\mu_0)=e^{\frac{1}{2\sigma^2}(\mu - \mu_0)^2 - \frac{1}{\sigma^2} (y-\mu_0)(\mu - \mu_0)}$ 
\end{small}
for normal error, while for Laplace error, 
\begin{small} 
h$_\sigma(y,\mu-\mu_0)=e^{\frac{1}{\sigma}(|\mu_0 - \mu| )}$.
\end{small} 
Further \small{$sup_x|\log(h_\sigma(y,\mu(x)-\mu_0(x)))|$ $\rightarrow0$} for all $y \in \Re$ as $||\mu-\mu_0||_\infty$, $\sigma^2$ go to 0 with $|| \mu - \mu_0 ||_{\infty}/\sigma^2 \to 0$. Under regularity conditions, $\frac{f_0}{f_{\mu,\sigma}}$ has a finite upper bound, hence we can use dominated convergence theorem subsequently. Observe, for a fixed $\sigma$, 
\begin{small}
$\int_{0}^{1} \Gamma_{\sigma}(y-\mu(x)) dx \ge  
\frac{1}{ sup_{x\in\left(0, 1\right)} h_\sigma\left(y,\mu\left(x\right)-\mu_{0}\left(x\right)\right)}\int_{0}^{1}\Gamma_{\sigma}(y-\mu_{0}(x))$dx.
\end{small}
As $||\mu-\mu_0||_\infty$, $\sigma^2$ go to 0 with $|| \mu - \mu_0 ||_{\infty}/\sigma^2 \to 0$, and applying dominated convergence theorem, $0<\int_{\Re} f_0(y)\log\frac{f_0(y)}{f_{\mu,\sigma}(y)} dy \le \lim_{\sigma \rightarrow 0}\lim_{||\mu-\mu_0||_\infty \rightarrow 0}\log(sup_x h_\sigma(y, \mu(x)-\mu_0(x))) \to 0$. For Cauchy errors, \small{$f_0(y)=\lim_{\sigma \to 0}\int_0^1\frac{1}{\pi\sigma}\frac{1}{\bigg(1+\frac{(y-\mu_0(x))^2}{\sigma^2}\bigg)}dx$} and \\ \small{$f_{\mu,\sigma}(y)=\int_0^1\frac{1}{\pi\sigma}\frac{1}{\bigg(1+\frac{(y-\mu(x))^2}{\sigma^2}\bigg)}dx=
\int_0^1\frac{1}{\pi\sigma}\frac{1}{\bigg(1+\frac{1}{\sigma^2}\left[(y-\mu_0(x))^2-2(y-\mu_0(x))(\mu(x)-\mu_0(x))+(\mu(x)-\mu_0(x))^2\right]\bigg)}dx$}. 
\begin{normalsize}
As $||\mu-\mu_0||_\infty$, $\sigma^2$ go to 0 with $|| \mu - \mu_0 ||_{\infty}/\sigma^2 \to 0$, and applying dominated convergence theorem, $\int_{\Re} f_0(y)\log\frac{f_0(y)}{f_{\mu,\sigma}(y)}dy \to 0$.
Thus taking appropriate limits, KL($f_{\mu,\sigma}, f_0$) goes to 0 under Gaussian, Laplace and Cauchy errors. Hence we can choose a suitably small $\epsilon_{1}$ and $\epsilon_{2}$ with $0<\epsilon_{1}<\epsilon_{2}$ such that \small{$\bigg\{||\mu-\mu_{0}||_\infty \le \epsilon_{1}, 0<\sigma \le \epsilon_{2}\bigg\} \Rightarrow KL(f_{\mu,\sigma}, f_0) \le \epsilon$}. Using the assumptions on the support of priors $\Pi^{*}$ and $\nu$, we have $\Pi(KL_\epsilon(f_0))>0$.
\end{normalsize}
\vskip12pt

\textbf{Proof of Theorem 2}:\\
\begin{normalsize}
Using the regularity conditions and Taylor's series expansion, we have for fixed y, $\mu$, $\sigma$,
\begin{eqnarray*}
\log\frac{f_0(y)}{f_{\mu,\sigma(y)}}&=&\sum_{k=1}^{n_0}\left(-1\right)^k\frac{\left\{(f_0(y)-1)^k-(f_{\mu,\sigma}(y)-1)^k\right\}}{k} + \delta^y_1\left(n_0\right)-\delta^y_2\left(n_0\right),\mbox{ for a fixed n}_0, 
\end{eqnarray*}
where $\delta^y_1\left(n_0\right)-\delta^y_2\left(n_0\right)$ is uniformly bounded in y. Using the identity a$^{n}$-b$^{n}$=(a-b)($\sum_{k=1}^n a^{n-k}b^{k-1}$), and denoting $g_0=f_0-1$ and $g_{\mu,\sigma}=f_{\mu,\sigma}-1$, we have, 
\begin{eqnarray*}
\int|\sum_{k=1}^{n_0}\left(-1\right)^{k}\frac{\left\{(f_0-1)^k-(f_{\mu,\sigma}-1)^k\right\}}{k} |dy \le
\sum_{k=1}^{n_0}\int|\frac{(-1)^k}{k}(f_0-f_{\mu,\sigma})(\sum_{l=1}^k g_{\mu,\sigma}^{k-l}g_0^{l-1})|dy\\
\le\sum_{k=1}^{n_{0}}\sup_y|\frac{(-1)^k}{k}(\sum_{l=1}^k g_{\mu,\sigma}^{k-l}g_0^{l-1})|\int|f_0(y)-f_{\mu,\sigma}(y)|dy=K(n_0)\int|f_0(y)-f_{\mu,\sigma}(y)|dy, 
\end{eqnarray*}
where K(n$_0$)=$\sum_{k=1}^{n_{0}}\sup_y|\frac{(-1)^k}{k}(\sum_{l=1}^k g_{\mu,\sigma}^{k-l}g_0^{l-1})|$ is a finite constant depending on n$_0$, for $\mu$ belonging to a finite L-1 ball around $\mu_0$, using the regularity conditions. Further, using similar methods as in the proof of theorem $3$, we can show that for $\epsilon_1<\epsilon_2<\epsilon_2^*$,
\begin{eqnarray}
\left\{\mu\in N_{\epsilon_1}(\mu_0), \sigma\in(\epsilon_2, \epsilon_2^*) \right\}\Rightarrow\int |f_0(y)-f_{\mu,\sigma}(y)|dy<\frac{\epsilon_{1}}{\epsilon_{2}}. \label{eq:L1}
\end{eqnarray}
Using inequality ($\ref{eq:L1}$), we have for $\mu\in N_{\epsilon_1}(\mu_0)$ and $\sigma\in(\epsilon_2, \epsilon_2^*)$,
\begin{eqnarray*}
&\int&|\sum_{k=1}^{n_0}\left(-1\right)^{k}\frac{\left\{(f_0-1)^k-(f_{\mu,\sigma}-1)^k\right\}}{k} |dy \le K(n_0)\frac{\epsilon_1}{\epsilon_2}. \mbox{ Also note, } KL(f_0,f_{\mu,\sigma}) \le \int f_0|\log\frac{f_0}{f_{\mu,\sigma}}|\\
&=& \int f_0|\sum_{k=1}^{n_0}\left(-1\right)^k\frac{\left\{(f_0-1)^k-(f_{\mu,\sigma}-1)^k\right\}}{k} + \delta^y_1\left(n_0\right)-\delta^y_2\left(n_0\right)| \le K(n_0)\frac{\epsilon_1}{\epsilon_2}+\Delta(n_0)=\epsilon,
\end{eqnarray*}
for a finite K$(n_0)$ and suitably small $\Delta(n_0)=sup_y|\delta^y_1(n_0)-\delta^y_2(n_0)|$ with $\epsilon_1$, $\epsilon_2$ depending on $\epsilon$. 
Under positive L-1 support by $\Pi^*$, the rest follows by similar arguments as in theorem 1.
\end{normalsize}
\vskip12pt

\begin{normalsize}
\textbf{Proof of Theorem 3}:
\begin{eqnarray*}
\mbox{Note that, } |f_{\mu,\sigma}(y)-f_0(y)| \le  |sup_{x \in (0,1)}\phi_{\sigma}(y-\mu(x))-f_0|, \mbox{ which is integrable } \forall (\mu,\sigma).
\end{eqnarray*}
Then using dominated convergence theorem, $\lim_{\sigma \to 0}\int |f_{\mu,\sigma}-f_0|=\int \lim_{\sigma \to 0} |f_{\mu,\sigma}-f_0|$ .
Now applying Fatou's Lemma and Fubini's Theorem successively, we have,
\begin{eqnarray*}
\int \lim_{\sigma\rightarrow0}|f_{\mu,\sigma}-f_0|dy &\le&  \int \lim_{\sigma\rightarrow0}\int_{0}^{1}|\phi_\sigma(y-\mu(x))-\phi_{\sigma}(y-\mu_{0}(x))|\mbox{dx dy}\\
&\le& \liminf_{\sigma\rightarrow0} \int\int_{0}^{1}|\phi_\sigma(y-\mu(x))-\phi_{\sigma}(y-\mu_{0}(x))|\mbox{dx dy}\mbox{ (Fatou's lemma) }\\
&=& \liminf_{\sigma\rightarrow0} \int_{0}^{1}\int|\phi_\sigma(y-\mu(x))-\phi_{\sigma}(y-\mu_{0}(x))|\mbox{dy dx}\mbox{ (Fubini's Theorem) }\\
&=& \liminf_{\sigma\rightarrow0}\bigg\{ \int_{x|\mu_{0}>\mu}\int|\phi_\sigma(y-\mu(x))-\phi_{\sigma}(y-\mu_{0}(x))|\mbox{dy dx}\\
&+& \int_{x|\mu_0<\mu}\int|\phi_{\sigma}(y-\mu_{0}(x))-\phi_\sigma(y-\mu(x))|\mbox{dy dx}\bigg\}.
\end{eqnarray*}
In the proof of lemma $1$ of Ghosal, Ghosh and Ramamoorthy ($1999$), it was shown that for fixed $\theta_1<\theta_2$, \small{$||\phi_\sigma(y-\theta_1)-\phi_\sigma(y-\theta_2)||<\frac{\theta_2-\theta_1}{\sigma}$}, which would imply 
\begin{eqnarray}
\int \lim_{\sigma\rightarrow0}|f_{\mu,\sigma}-f_0|dy &\le& \liminf_{\sigma\rightarrow0}\bigg\{ \int_{x|\mu_{0}>\mu} \frac{\mu_0(x)-\mu(x)}{\sigma}dx + \int_{x|\mu_0<\mu} \frac{\mu(x)-\mu_0(x)}{\sigma}dx \bigg\} \nonumber \\
&=& \liminf_{\sigma\rightarrow0} \int_0^1 \frac{|\mu(x)-\mu_0(x)|}{\sigma} dx \le \lim_{\sigma\rightarrow0} \int_0^1 \frac{|\mu(x)-\mu_0(x)|}{\sigma} dx. \label{eq:CRUX}
\end{eqnarray}
As $||\mu-\mu_0||_\infty$, $\sigma^2$ go to 0 with $|| \mu - \mu_0 ||_{\infty}/\sigma^2 \to 0$, the above limit exists and goes to 0. Given that the limit exists and goes to $0$, we can now choose sufficiently small $\epsilon_{1},\epsilon_{2},\epsilon_{2}^{*}$ with $0<\epsilon_{1}<\epsilon_{2}<\epsilon_{2}^{*}$ such that for $\int_0^1|\mu_0(x)-\mu(x)|dx<\epsilon_1 \equiv \mu \in N_{\epsilon_1}(\mu_0)$ and $\sigma \in(\epsilon_{2},\epsilon_2^*)$, we would have $\int |f_{\mu,\sigma}-f_{0}|dy<\frac{\epsilon_{1}}{\epsilon_{2}}$, using ($\ref{eq:CRUX}$) . 
\end{normalsize}
\vskip12pt

\begin{normalsize}
\textbf{Proof of Theorem 4}:\\
Our proof is based on theorem $2$ of Ghosal, Ghosh and Ramamoorthi ($1999$) who gave a set of alternate sufficient conditions for almost sure convergence of the posterior of strong neighborhoods. Their result involves conditions on the size of the parameter space in terms of L-1 metric entropy. Before proceeding, let us review L-1 metric entropy and theorem $2$ of Ghosal, Ghosh and Ramamoorthi ($1999$). 

\begin{footnotesize}\textbf{DEFINITION 1}.\end{footnotesize} For $\mathcal{G}\subset \mathcal{F}$ and $\delta>0$, L-1 metric entropy $J(\delta,\mathcal{G})$ is defined as the minimum of $\log(k:\mathcal{G} \subset \cup_{i=1}^k \left\{f:\int|f-f_i|dy<\delta, \mbox{ } f_1,f_2,\ldots,f_k \in \mathcal{F}\right\} )$.

\begin{theorem}(Ghosal, Ghosh and Ramamoorthi) Let $\Pi$ be a prior on $\mathcal{F}$. Suppose $f_0 \in \mathcal{F}$ is in the Kullback-Leibler support of $\Pi$ and let U=$\left\{f:\int|f-f_0|dy<\epsilon\right\}$. If there is a $\delta<\epsilon/4$, $c_1,\mbox{ }c_2>0,\mbox{ }\beta<\epsilon^2/8$ and $\mathcal{F}_n\subset\mathcal{F}$ such that for all large n:\\
(1) $\Pi(\mathcal{F}_n^c)<c_1\exp(-nc_2)$, and, \\
(2) The L-1 metric entropy, $J(\delta,\mathcal{F}_n)<n\beta$, \\
then $\Pi(U|Y_1,Y_2,\ldots,Y_n)\rightarrow1$ a.s. $P_{f_0}$.
\end{theorem}
The constants $\delta$, $c_1,c_2,\mbox{ }\beta$ and $\mathcal{F}_n$ are allowed to depend on $\epsilon$.

Let the parameter space for $(\mu,\sigma)$ be denoted as $\mathcal{H}$. Consider the subsets of the parameter space $\mathcal{H}_n=\mathcal{H}_{1n}\otimes\mathcal{H}_{2n}$, where $\mathcal{H}_{1n}=\left\{\mu: ||\mu||_\infty< M_n, ||\mu'||_\infty < M_n\right\}$ and $\mathcal{H}_{2n}= [L_n,\infty)$, with $L_n \to$ 0 such that $\nu(\sigma \in (0,L_n))<d_1\exp(-d_2n)$, $d_1,d_2>0$ and $M_n=O(n^{1/2})$. The regularity conditions on the GP guarantees existence of the first derivative $\mu'$ with probability 1. Using lemma 4 of Choi and Schervish (2004) who showed an upper bound on sup-norm metric entropy of $\mathcal{H}_{1n}$, we have the upper bound on L1 metric entropy as  
\begin{eqnarray}
J(\delta,\mathcal{H}_{1n})<K_1M_n/\delta.   \label{eq:stA}
\end{eqnarray}
This implies there are K$^*=\exp(K_1 M_n/\delta)$ elements $\mu_1,\mu_2,\ldots,\mu_{K^*}$ such that 
\begin{eqnarray}
\mathcal{H}_{1n} \subset \cup_{j=1}^{K^*} \left\{\mu:\int_0^1|\mu-\mu_j| dx< \delta\right\}. \label{eq:stB}
\end{eqnarray}
Let us consider the sieve $\mathcal{F}_n=\left\{f_{\mu,\sigma} \in \mathcal{F}: (\mu,\sigma) \in \mathcal{H}_n\right\}$. Further, let us consider densities f$_{i,n}=f_{\mu_i,L_n}\in \mathcal{F}_n$ defined as in section 2, where the $\mu_i, i=1,\ldots, K^*$ correspond to the ones just defined to cover $\mathcal{H}_{1n}$. Using similar techniques as in lemma 1 of Ghosal, Ghosh and Ramamoorthi (1999) and theorem 3, it can be shown that for $f_{\mu,\sigma} \in \mathcal{F}_n$,
\begin{eqnarray}
\int |f_{\mu,\sigma}-f_{i,n}| dy &\leq& \int_0^1 \frac{|\mu-\mu_i|}{\sqrt{2\pi}\sigma} dx + \int_0^1 \frac{|\mu-\mu_i|}{\sqrt{2\pi}L_n} dx  \leq \frac{\sqrt{2}}{\sqrt{\pi}L_n} \int_0^1 |\mu-\mu_i| dx \leq \frac{\delta}{L_n}. \label{eq:stC}
 \end{eqnarray}
This clearly implies 
\begin{eqnarray}
J(\delta/L_n,\mathcal{F}_n) &=& J(\delta,\mathcal{H}_{1n}) <K_1M_n/\delta  \Rightarrow J(\delta,\mathcal{F}_n) \leq K_1M_nL_n/\delta  < n\beta, \label{eq:stD}
\end{eqnarray}
where we can choose $\delta<\epsilon/4$ such that $\beta<\epsilon^2/8$. Thus the second condition in theorem 6 is satisfied.
Also note that the prior probability of $\mathcal{F}_n$ can be calculated in terms of $\Pi^*$ and $\nu$. Using the regularity conditions on the GP and $\nu$, and similar reasoning as in lemma 5 of Choi and Schervish (2004),
\begin{eqnarray}
\Pi(\mathcal{F}_n^c) = (\Pi^*\otimes\nu)(\mathcal{H}_n^c) \leq c_1\exp(-c_2n), c_1,c_2>0. \label{eq:stE}
\end{eqnarray}
Thus the first condition in theorem 6 is satisfied.
\end{normalsize}

\begin{figure}
\centering
\mbox{\includegraphics[height=7.5 in, width= 1.0 \textwidth]{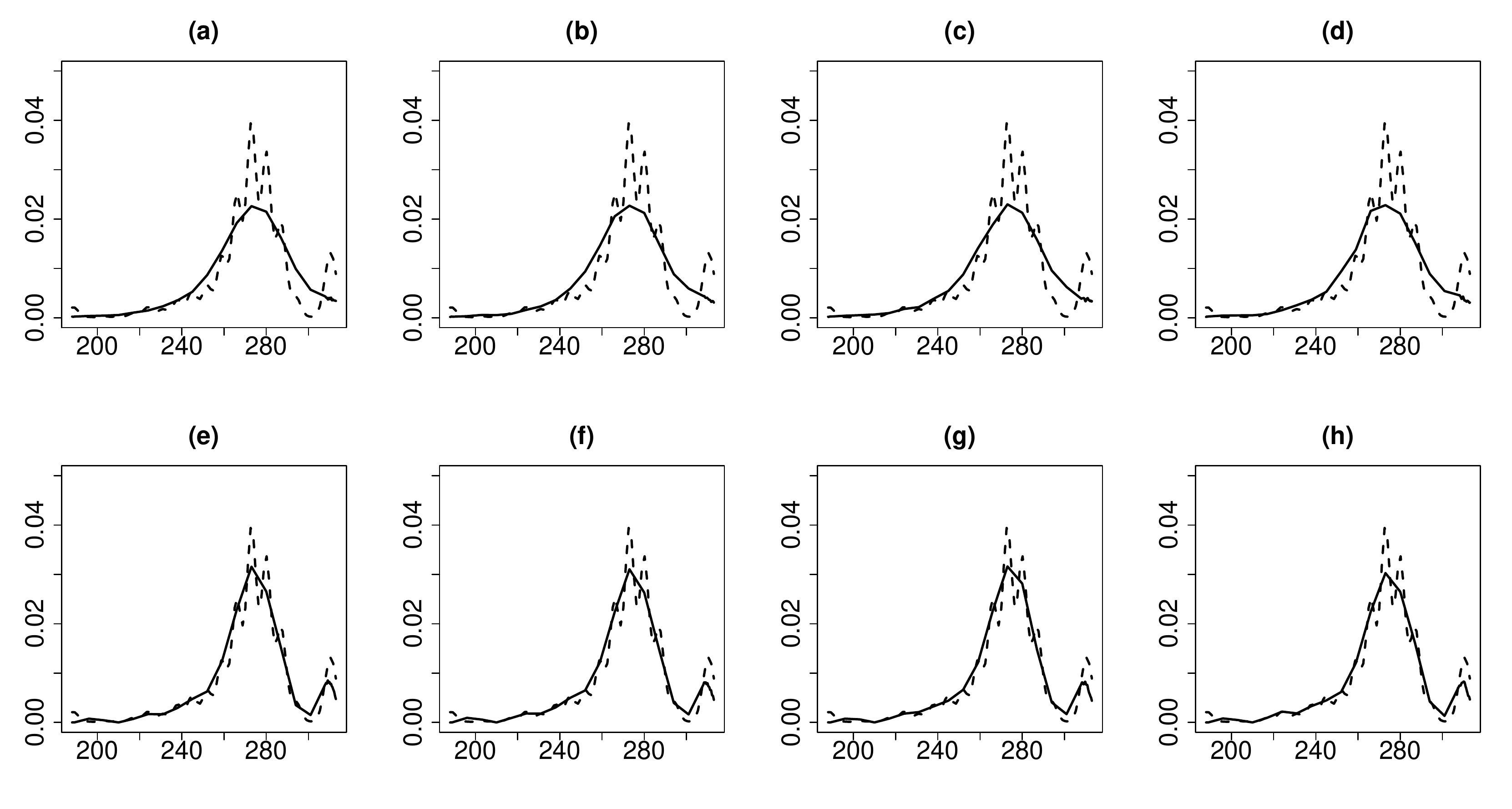}}
\caption{\small{Prior realizations from the GPT for gestational age at delivery (solid lines) along with frequentist kernel density estimate (dotted lines). The rows correspond to $\phi_1$=(0.01, 0.1); the columns correspond to $\phi_2$=(0.1,1,25,100).}}
\end{figure}

\begin{figure}
\centering
		\mbox{\includegraphics[height=4.2in, width=0.45\textwidth]{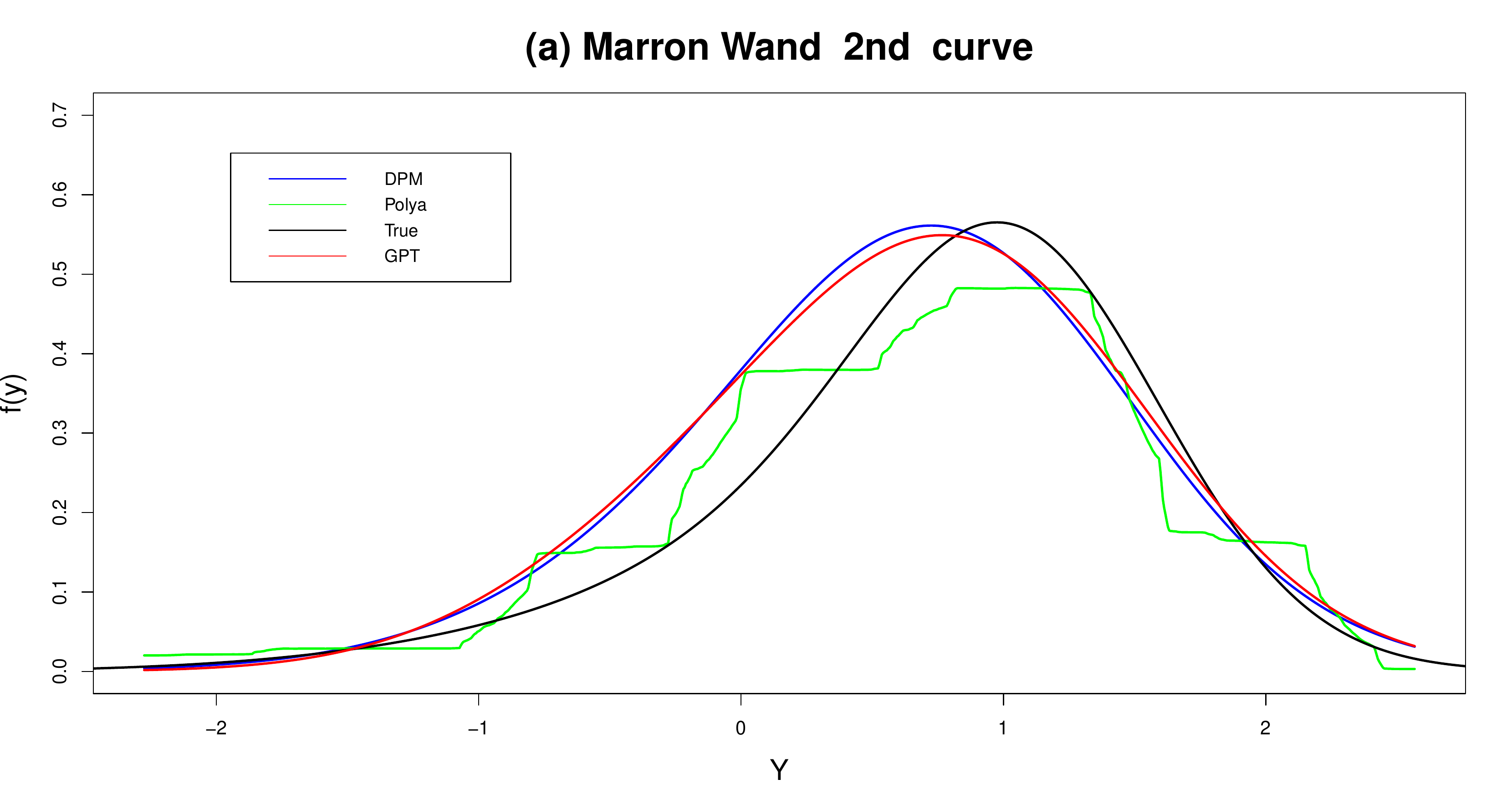}} \quad \mbox{\includegraphics[height=4.2in, width=0.45\textwidth]{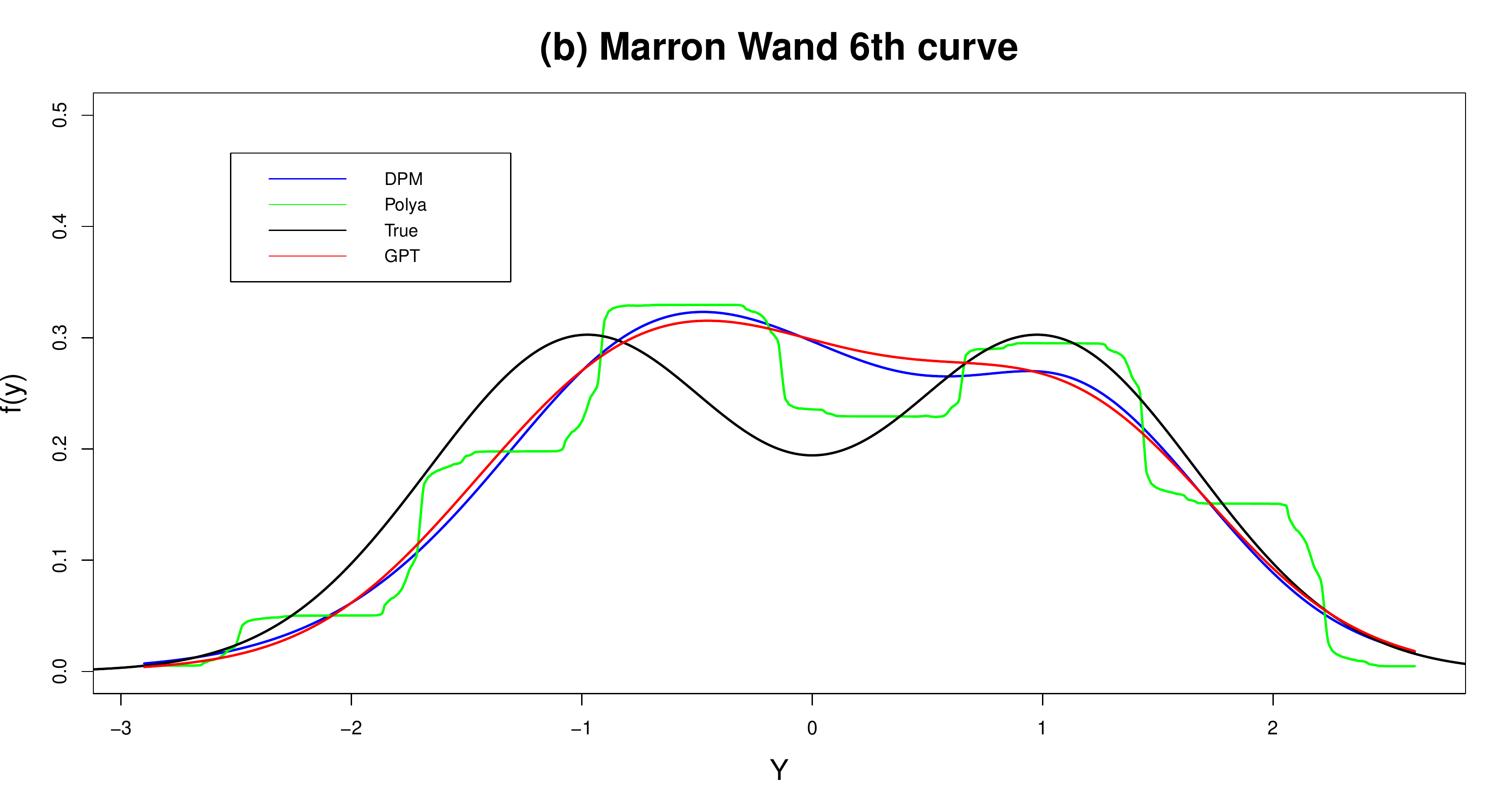}}\\
		\mbox{\includegraphics[height=4.2in, width=0.45\textwidth]{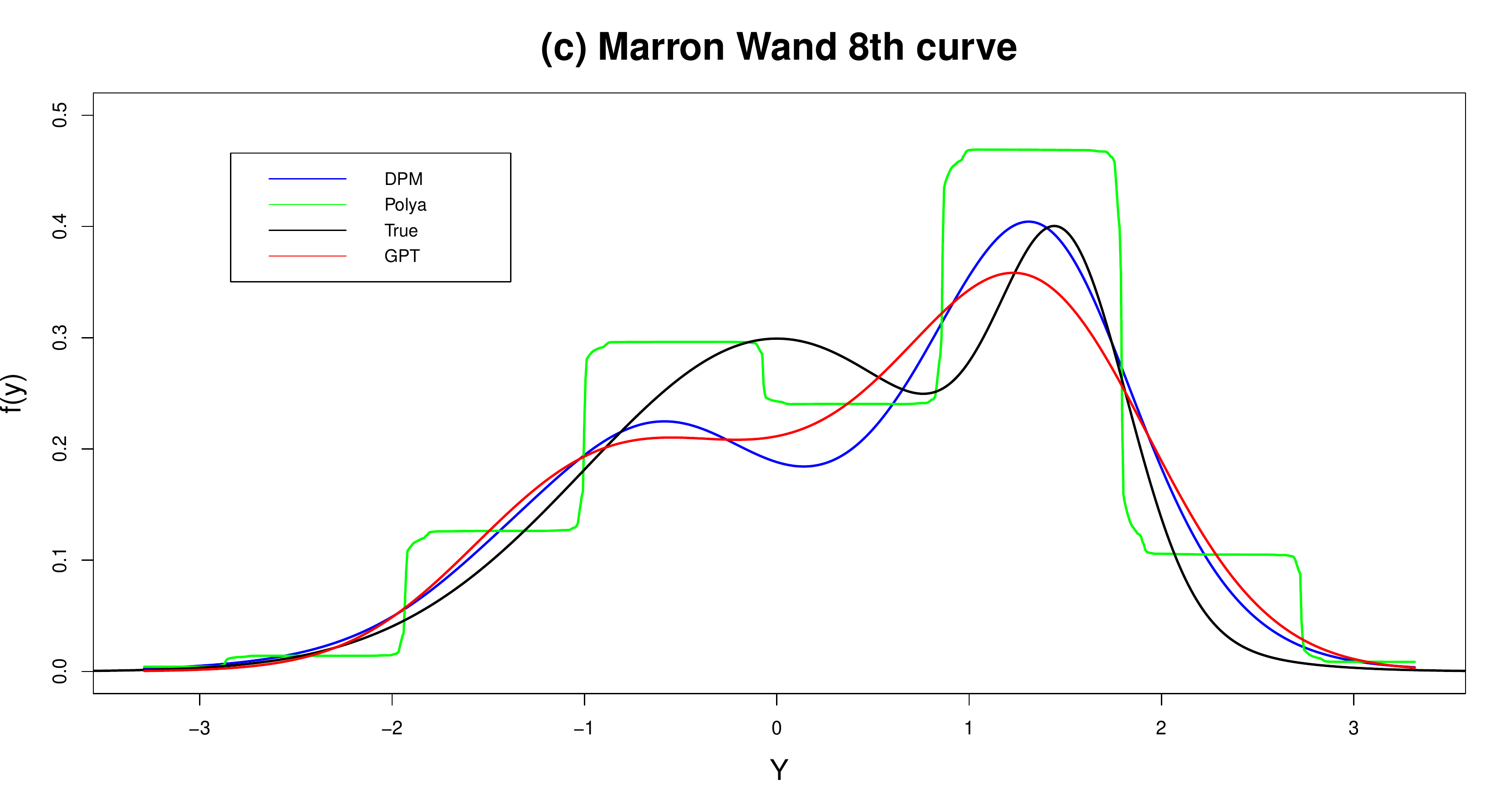}} \quad 	\mbox{\includegraphics[height=4.2in, width=0.45\textwidth]{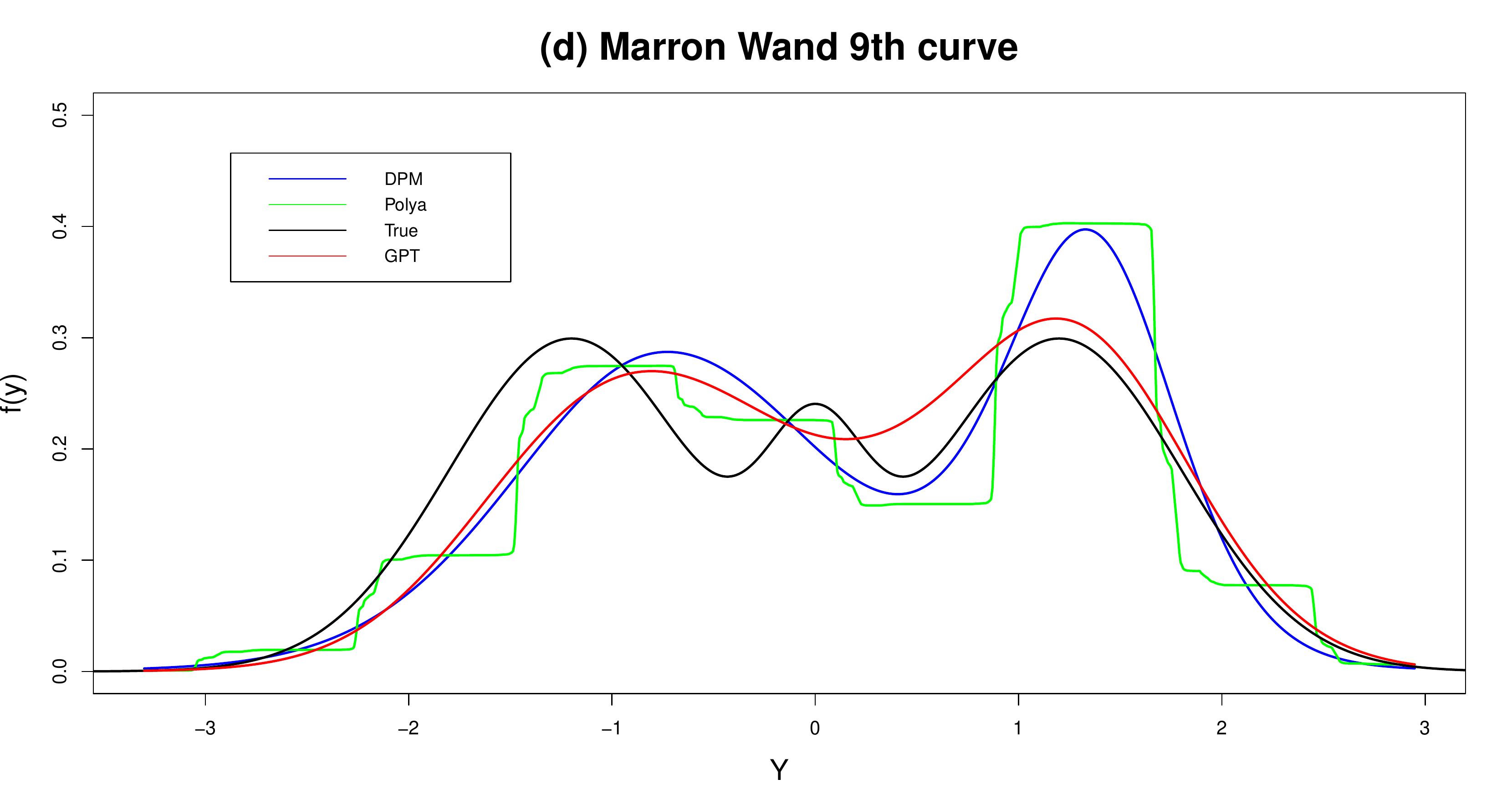}}
			\caption{\small{Marron-Wand curves - density estimates for GPT, DPM and Polya tree mixtures.}}
\end{figure}
\begin{figure}
\centering
		\mbox{\includegraphics[height=3.8in, width=3.0in]{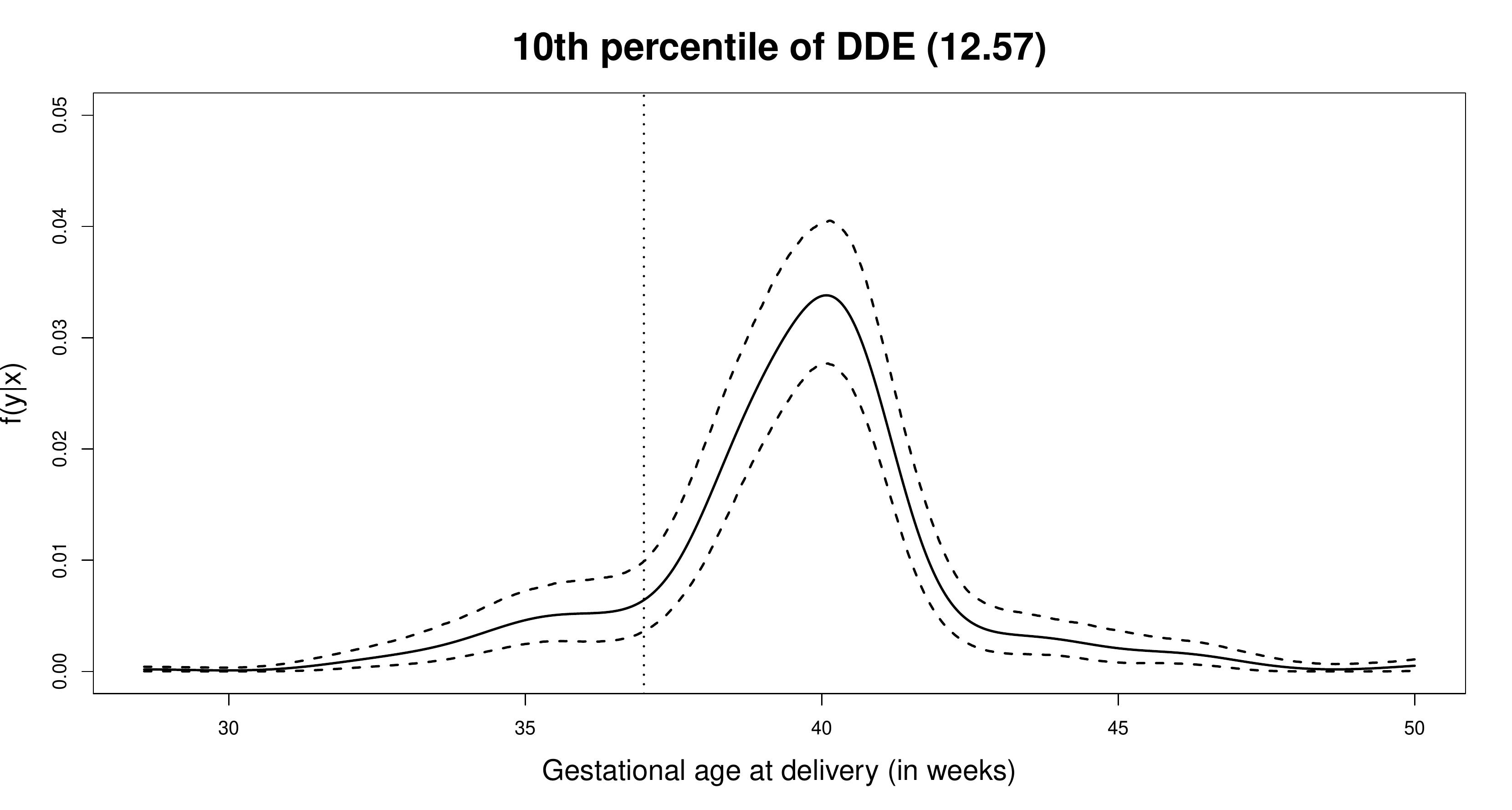}} \quad \mbox{\includegraphics[height=3.8in, width=3.0in]{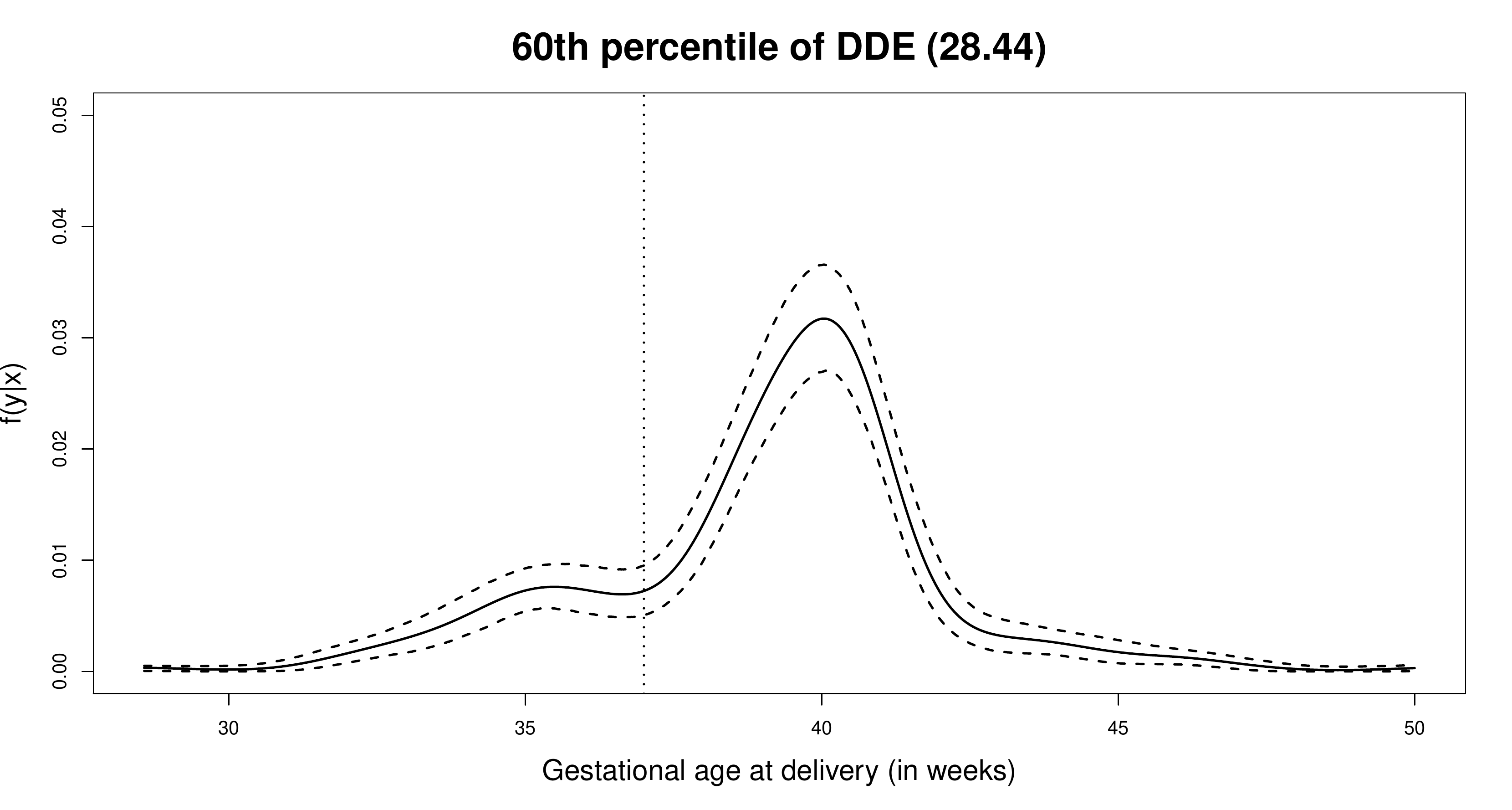}}\\
		\mbox{\includegraphics[height=3.8in, width=3.0in]{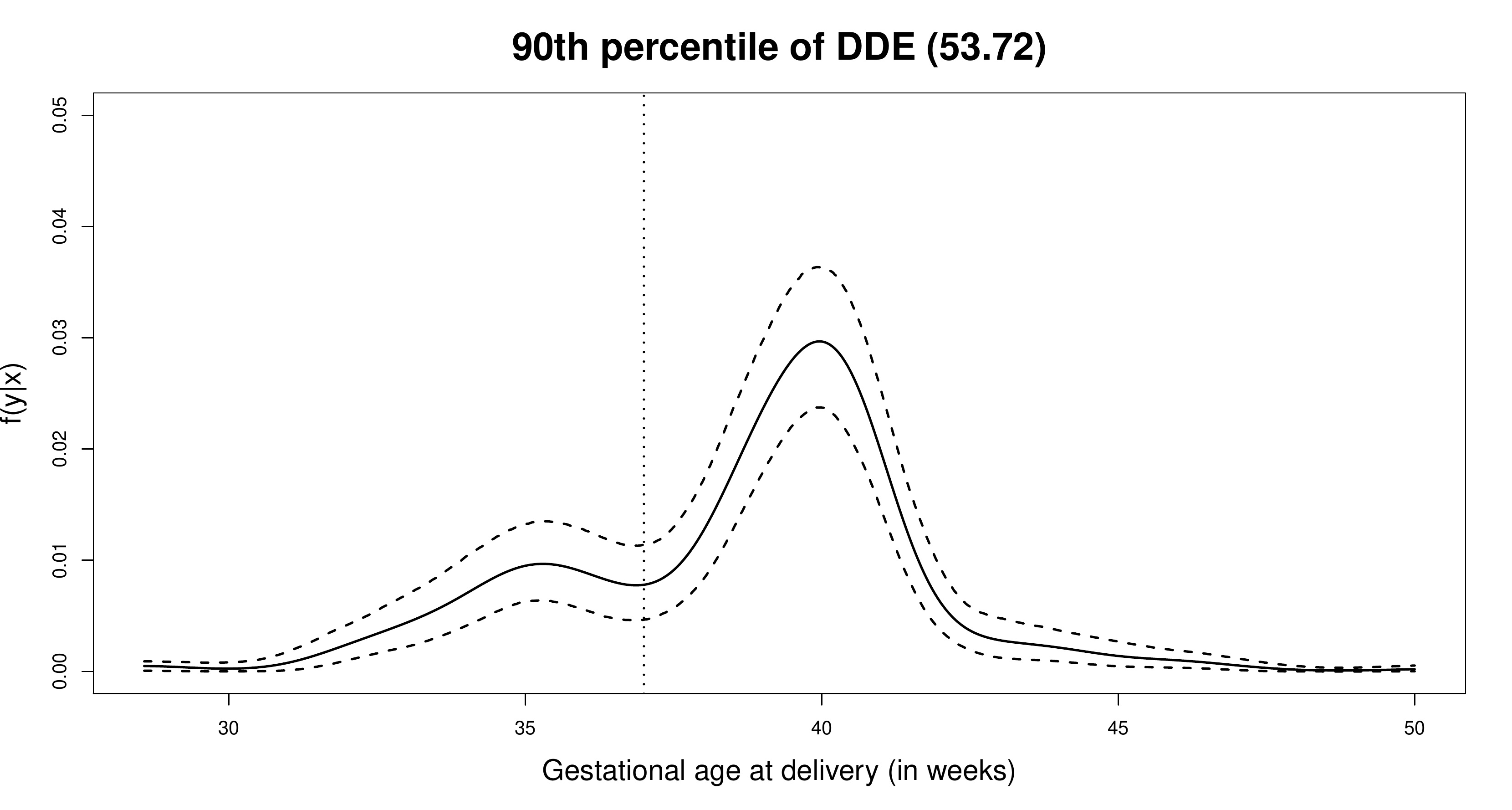}} \quad 	\mbox{\includegraphics[height=3.8in, width=3.0in]{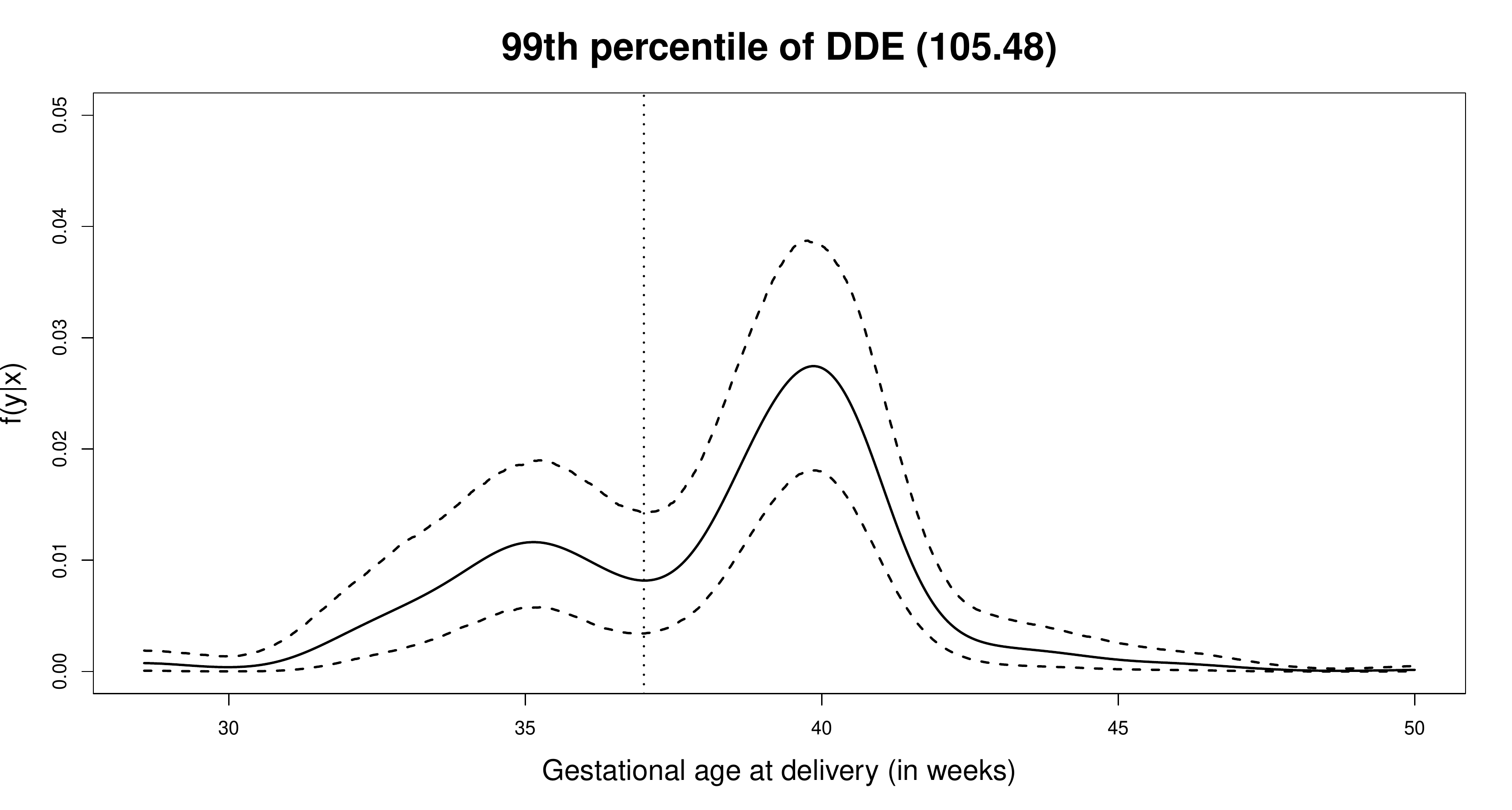}}\\
			\caption{\small{GPT conditional density estimates and 90\% credible intervals for 10th, 60th, 90th, 99th DDE quantiles.
			                Vertical dashed line for cut-off at 37 weeks.}}
\end{figure}
\begin{figure}
\centering
		\mbox{\includegraphics[height=3.8in, width=3.0in]{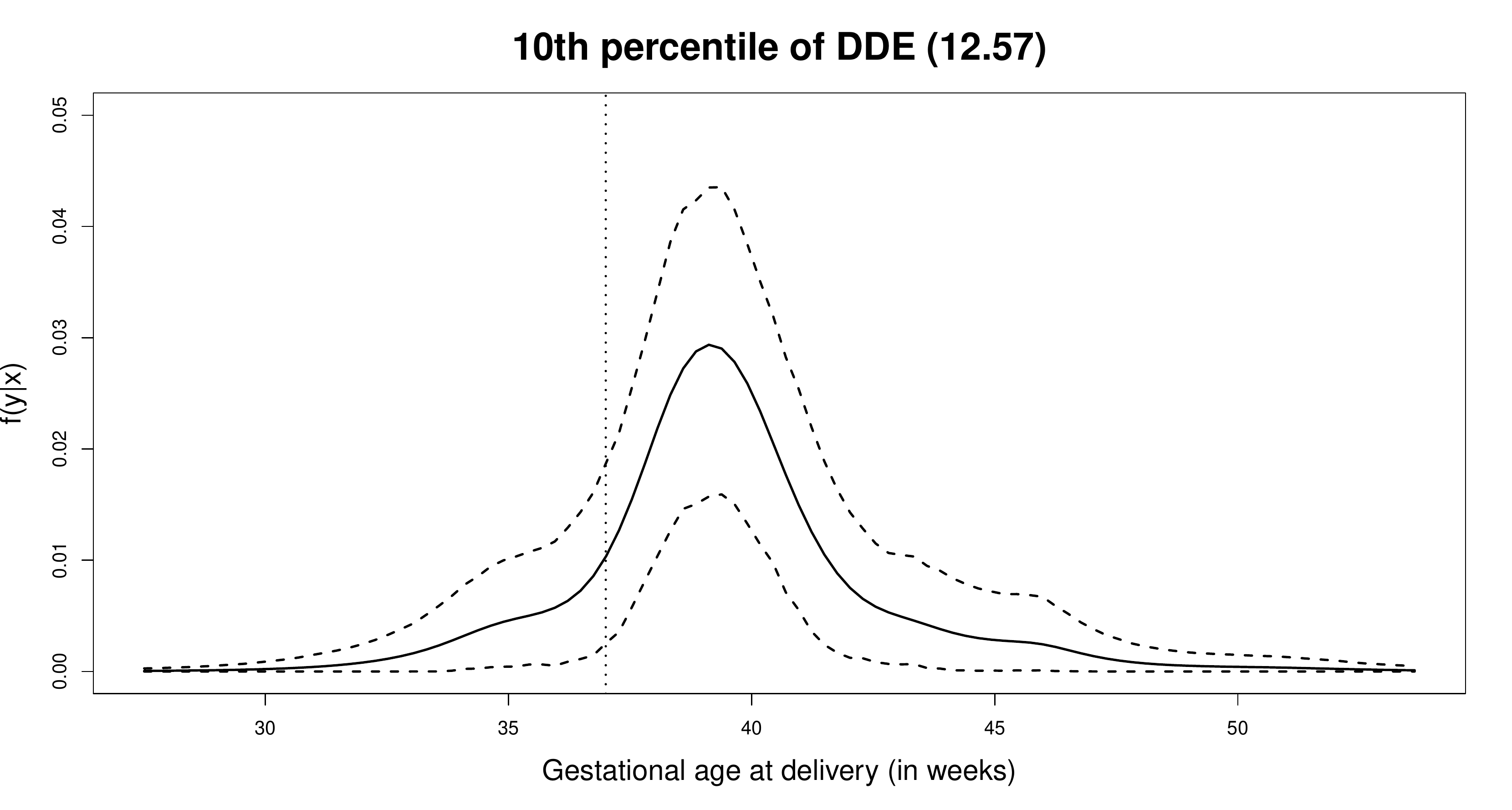} \quad \includegraphics[height=3.8in, width=3.0in]{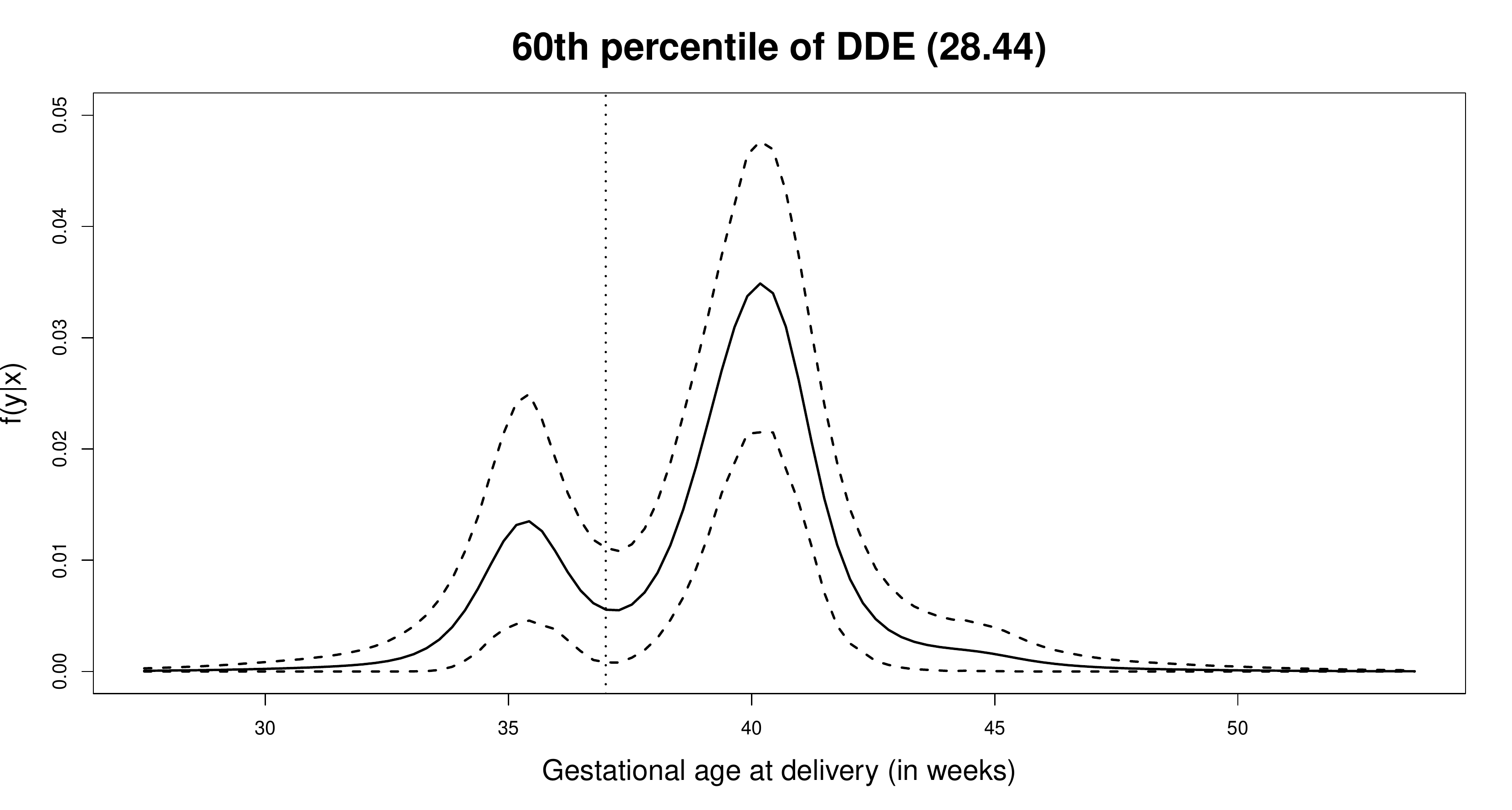}}\\
		\mbox{\includegraphics[height=3.8in, width=3.0in]{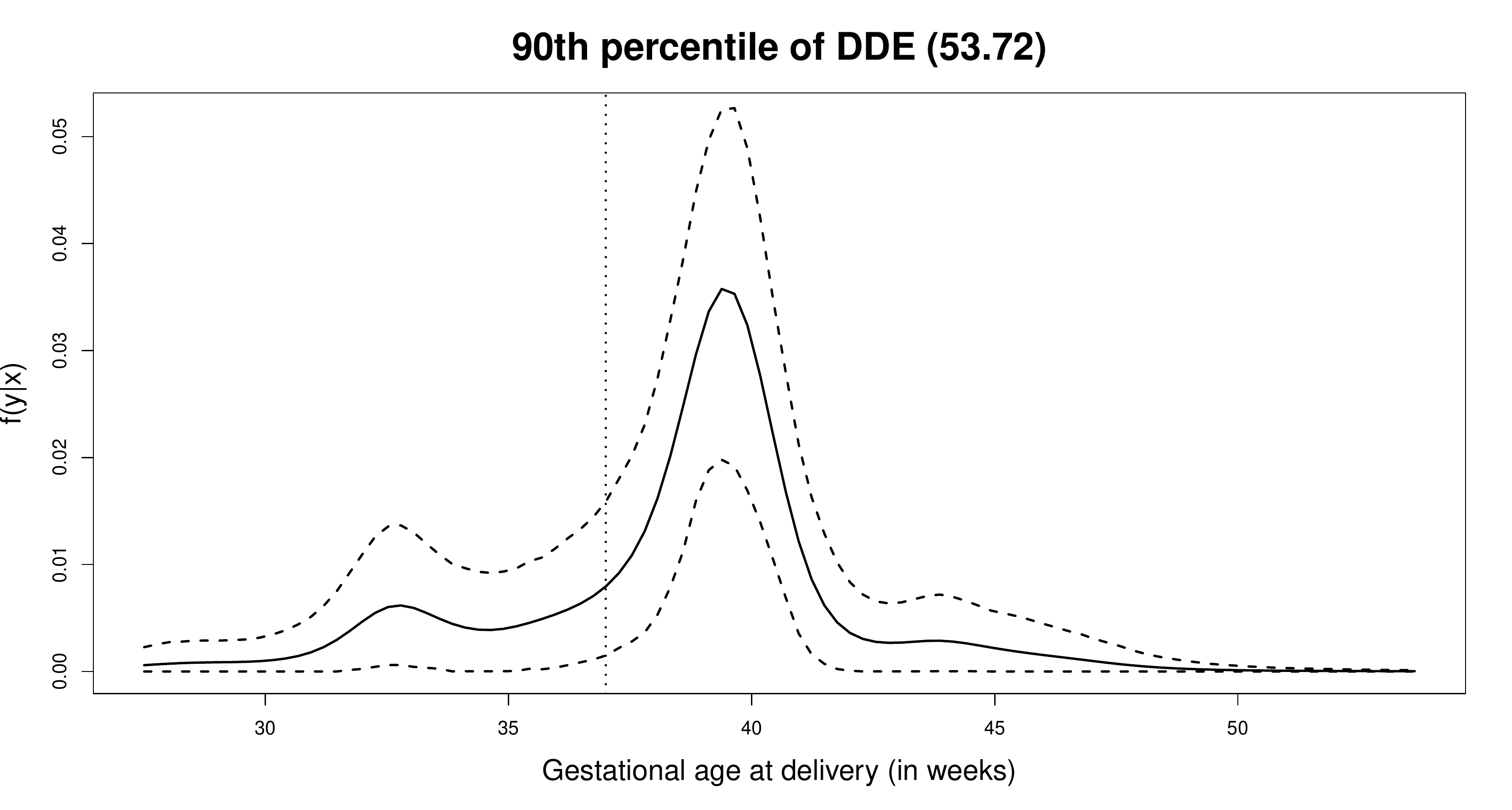} \quad \includegraphics[height=3.8in, width=3.0in]{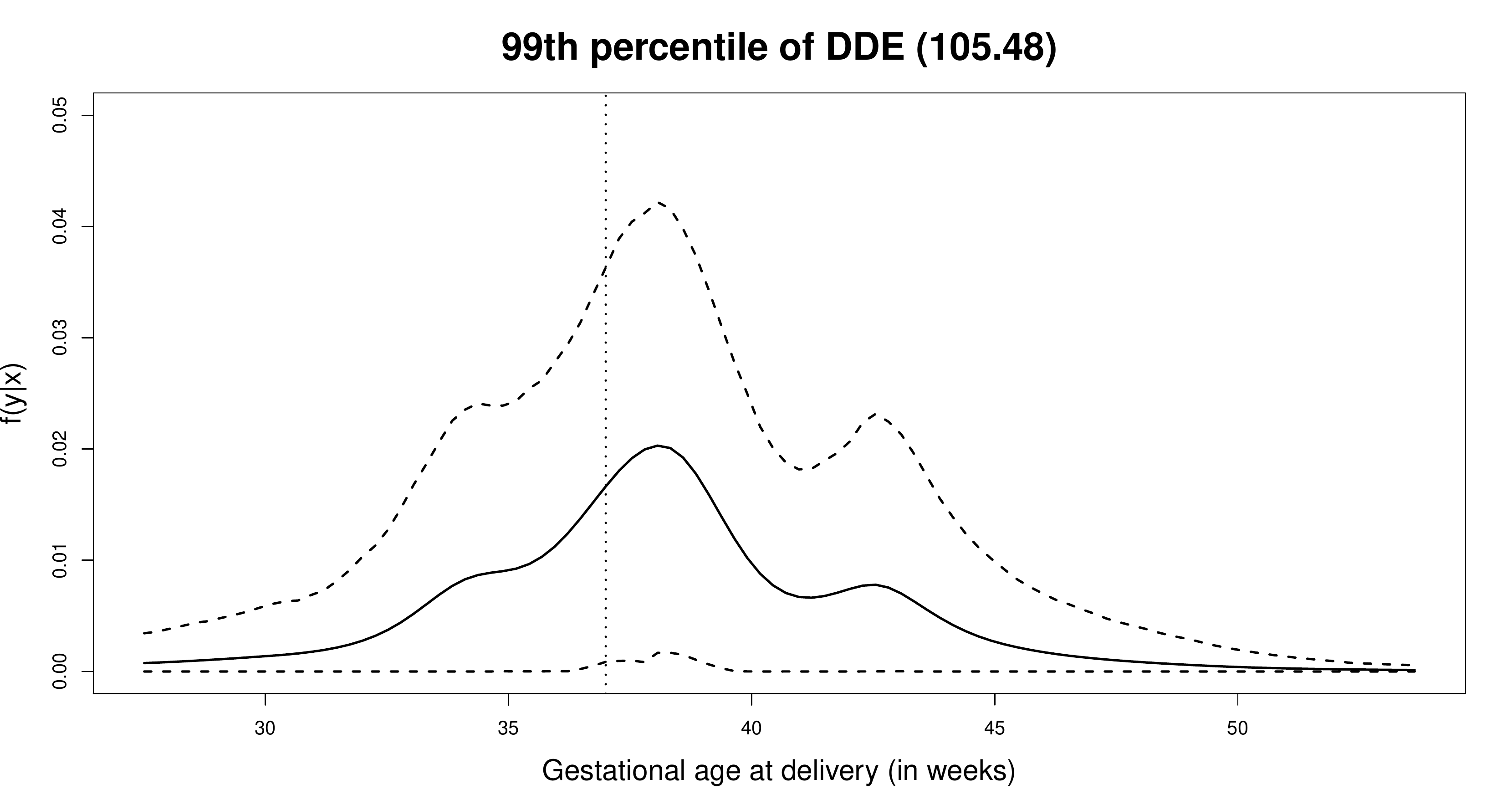}}\\
			\caption{\small{DPM conditional density estimates and 90\% credible intervals for 10th, 60th, 90th, 99th DDE quantiles.
			                Vertical dashed line for cut-off at 37 weeks.}}
\end{figure}
\begin{figure}
\centering
		\mbox{\includegraphics[height=3.8in, width=3.0in]{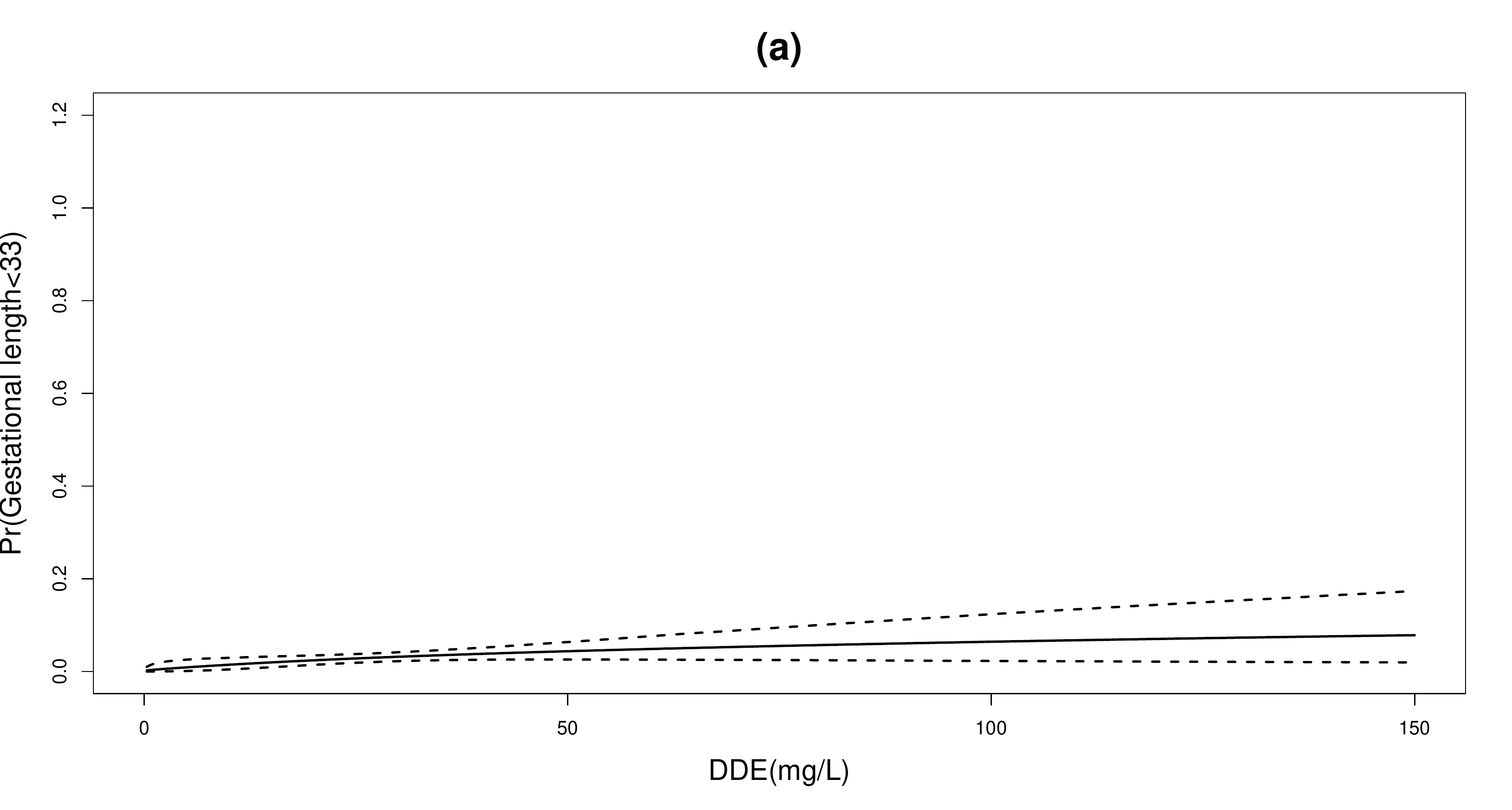}} \quad \mbox{\includegraphics[height=3.8in, width=3.0in]{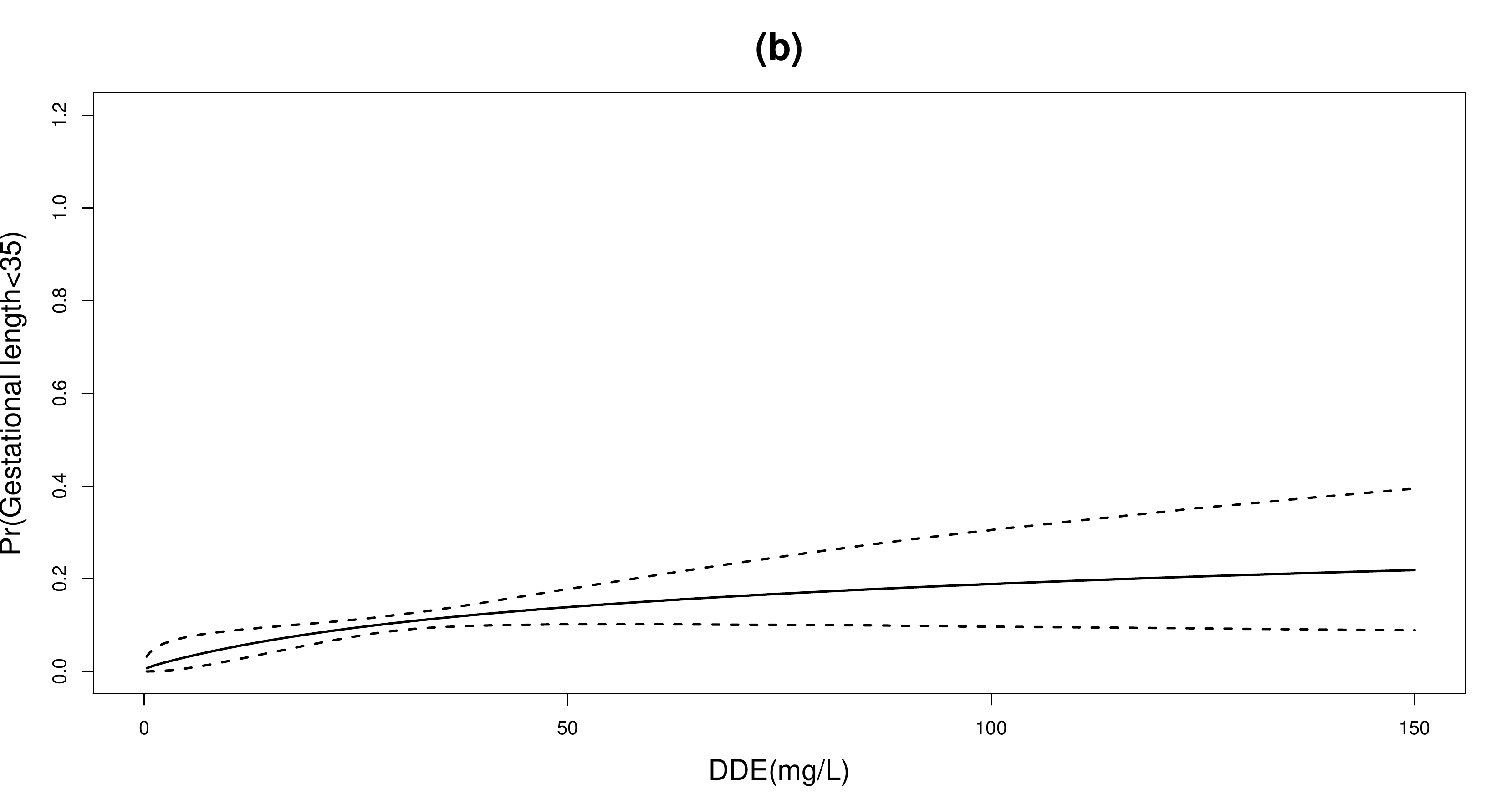}}\\
		\mbox{\includegraphics[height=3.8in, width=3.0in]{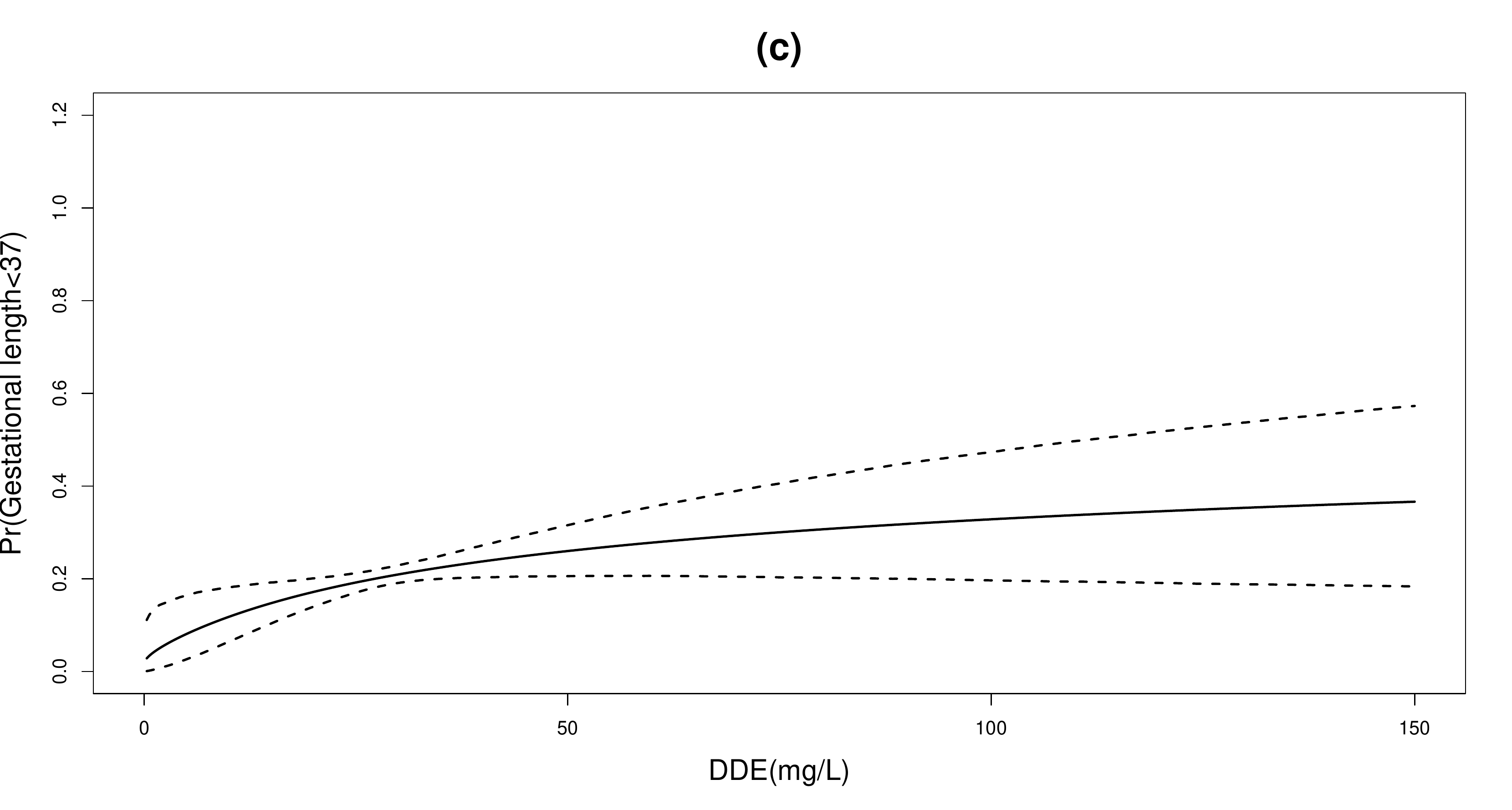}} \quad 	\mbox{\includegraphics[height=3.8in, width=3.0in]{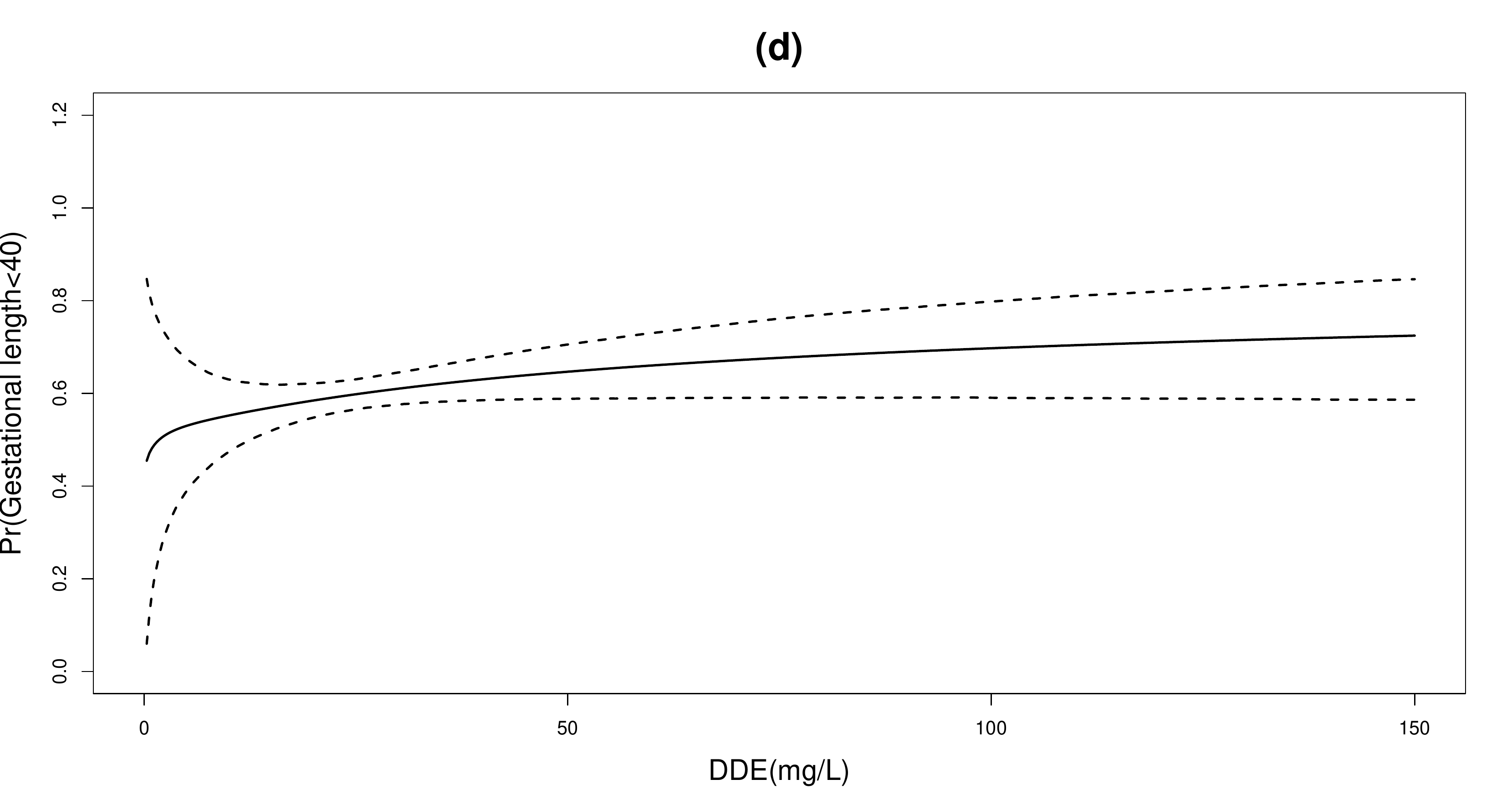}}\\
			\caption{\small{Estimated probability that gestational age at delivery is less than T weeks versus DDE dose, for
(a) T = 33, (b) T = 35, (c) T = 37, (d) T = 40. Solid lines are posterior means and dashed lines are pointwise 90\%
credible intervals.}}
\end{figure}

\newpage 

\linespread{0.8}
\begin{small}

\end{small}

\end{document}